\documentclass[a4paper,12pt]{amsart}

\usepackage{amsmath,amsthm}

\usepackage[pagewise]{lineno}

\usepackage{amsrefs,amssymb,amscd,xcolor}
\usepackage{tikz, subfigure}
\usepackage{mathrsfs}
\usepackage{mathtools}
\usepackage{dsfont}
\usepackage{upgreek}

\usepackage{a4wide}
\usepackage[normalem]{ulem}
\usepackage{MnSymbol}

\newtheorem{prop}{Proposition}[section]

\newtheorem{theorem}{Theorem}[section]
\newtheorem*{theorema}{Theorem A}

\newtheorem*{theoremb}{Theorem B}
\newtheorem*{theoremc}{Theorem C}
\newtheorem{lemma}{Lemma}[section]

\newtheorem{example}{Example}[section]
\newtheorem{remark}{Remark}[section]
\theoremstyle{definition}
\newtheorem{defn}{Definition}[section]

\newcommand{\R}{\mathbb{R}}

\newcommand{\N}{\mathbb{N}}

\newcommand{\E}{\mathbb{E}}
\newcommand{\Z}{\mathbb{Z}}
\newcommand{\I}{{\mathcal I}}

\newcommand{\dom}{\text{dom}}

\def\L{\mathcal{L}}

\newcommand{\A}{\mathcal A}

\newcommand{\B}{\mathcal B}

\newcommand{\La}{\mathcal L}
\newcommand{\M}{\mathcal M}
\newcommand{\e}{\mathcal E}
\newcommand{\Zb}{\mathcal Z}

\renewcommand{\top}{{\textsf{top}}}

\renewcommand{\phi}{{\varphi}}
\renewcommand{\epsilon}{{\varepsilon}}
\newcommand{\eps}{{\varepsilon}}
\def\mybf #1{\textbf{\textit{#1}}}

\author{Mirmukhsin Makhmudov}
\address[M.~Makhmudov]{Mathematical Institute\\
 Leiden University\\
 Postbus 9512, 2300 RA Leiden, The Netherlands}
\email{m.makhmudov@math.leidenuniv.nl}
\author{Evgeny Verbitskiy}
\address[E.~Verbitskiy]{Mathematical Institute\\
  Leiden University\\
   Postbus 9512, 2300 RA Leiden, The Netherlands}
 \address{Bernoulli Institute\\
  University of Groningen\\
   PO Box 407, 9700 AK Groningen, The Netherlands}
\email{evgeny@math.leidenuniv.nl}
\author{Qian Xiao}
\address[Q.~Xiao]{Mathematical Institute\\
  Leiden University\\
  Postbus 9512, 2300 RA Leiden, The Netherlands}
  \address{School of Mathematics\\
  South China University of Technology\\
  Guangzhou 510641, People's Republic of China}
\email{xiaoqianmath@163.com}

\begin{document}

\title{Multifractal Formalism from Large Deviations}
\date{\today}

\begin{abstract}	
It has often been observed that the Multifractal Formalism
and the Large Deviation Principles are intimately related. In fact, Multifractal Formalism 
was heuristically derived using the Large Deviations ideas. 
In numerous
examples in which the multifractal results have been rigorously established, 
the corresponding Large Deviation results are valid as well. 
Moreover, the proofs of multifractal and large deviations are remarkably similar.
The natural question then is whether under which conditions multifractal formalism 
can be deduced from the corresponding large deviations results.
More specifically, 
given a sequence of random variables $\{ {X_n} \}_{n\in\N}$, satisfying
a Large Deviation Principle, what can be said about the  multifractal nature
of the level sets $K_\alpha=\{\omega:  \lim_{n} \frac{X_n(\omega)}{n}=\alpha\}$. 
Under some technical assumptions, we establish the upper and lower bounds for multifractal spectra
in terms of the large deviation
rate functions, and show that many known results of multifractal formalism
are covered by our setup.

	
\end{abstract}

\maketitle

\section{Introduction}
The concept of multifractal formalism (MF) was proposed in 1980's  by Parisi, Frisch, Hentchel  \cites{Benzi-etc1984,Frisch-Parisi1983}
in the context of the study of turbulence.
It has been suggested  that natural local quantities have multiple fluctuation scales.
These possible fluctuations can be described using a \emph{singluarity spectrum}. The multifractal
formalism states that the singularity spectrum is dual, in the sense of Legendre transform, to some
integral ‘free energy'-type function of the system. The basic idea behind multifractal formalism 
is underpinned (or motivated) by an assumption that a certain probabilistic
the Large Deviation Principle (LDP) holds for the system.

The first rigorous mathematical results \cites{Olsen1995,Olsen1996,Barreira-Pesin-Schmeling1997} have been obtained in the early and mid-1990s (Riedli, Falconer, Olsen, Pesin), and since the late 1990's there was an explosion of research in this area \cites{Olsen2002,Olivier,Chazottes,Olivier1,Olivier2,Olivier3,Fan1996,Fan3,Fan1,Fan4}.
Multifractal formalism has been rigorously established for very large classes of dynamical systems
and a plethora of local observables. However, a somewhat curious phenomenon has occurred.
Though the original motivation and first rigorous results relied rather heavily on probabilistic 
methods of Large Deviations,
in subsequent works the link to Large Deviations has become somewhat less pronounced.
In fact, we are not aware of a single rigorous multifractal result without an accompanying 
Large Deviations result. 
However, it is also clear that more assumptions are required for the 
validity of multifractal formalism than for the validity of Large Deviations Principles.
For example, for multifractal analysis, the phase space must be a metric space, while this
is not really a requirement for Large Deviations.

This brings us to the natural question: assuming that the local quantity of interest
is a pointwise limit of quantities, whose probabilistic behaviour is 
governed by Large Deviations, what could be said about the corresponding singularity spectrum and multifractal formalism?

Let us compare two old results, first LDP for Bernoulli random variables, and second, what we would now
call a multifractal result for frequencies of digits in binary expansions.
However, this result precedes the multifractal idea by nearly 50 years.

\begin{theorem}[LDP for coin-tossing] Suppose $\left(X_k\right)$ are i.i.d. Bernoulli random variables with $\mathbb{P}\left(X_1=0\right)=\mathbb{P}\left(X_1=1\right)=$ $\frac{1}{2}$. Set $S_n=\sum_{i=1}^n X_i$, then $\{S_n/n\}$ satisfies LDP with the rate function $I(p)=\log2-H(p)$, where $H(p)=-p\log p-(1-p)\log (1-p)$ if $p\in[0,1]$, otherwise $H(p)=-\infty$. Namely,
$$
\lim_{n\to\infty} \frac{1}{n}\log {\mathbb{P}\Big(\frac{S_n}{n}\in(a,b)\Big)}=-\inf_{p\in(a,b)}I(p), \;\; \text{ for all } -\infty\leq a\leq b\leq\infty.
$$

\end{theorem}

The above theorem is one of the first Large Deviation results, and it is useful to compare it with possibly the first multifractal result by Besicovitch \cite{Besicovitch1935} and Eggleston \cite{Eggleston1949} in the 1930s and the 1940s. This multifractal result was, first, proved by Besicovitch \cite{Besicovitch1935} for $N=2$ (binary shift), and later, extended by Eggleston \cite{Eggleston1949} for general $N$ ($N$ is the cardinality of the alphabet).

\begin{theorem}[Besicovitch-Eggleston] Let $\Omega=\{0,1\}^{\Z_+}$, then
$$
 \operatorname{dim}_H\left\{\omega\in\Omega :\,
 \lim_n \frac 1n \sum_{j=0}^{n-1}\omega_j=p\right\}=
 \frac{H(p)}{\log 2}, \;\; \forall p\in[0,1],
 $$
 where $\dim_H$ denotes the Hausdorff dimension.
	
\end{theorem}
It is important to stress that not only rate functions in the above theorems are related, but in fact, the proofs are also very similar.

Next to the Hausdorff dimension, there is another, but in some sense more dynamical set characteristic, known as \textit{topological entropy}, which is used frequently to analyse non-compact non-invariant sets. However, in the setting of symbolic dynamics $(\Omega:=(0,...,l-1)^{\Z_+},\sigma)$, one has $h_{\top}(\sigma,Z)=\log l\cdot \dim_{H}(Z)$ for all $Z\subset \Omega$, in particular, $h_{\top}(\sigma,\Omega)=\log l$, where $\sigma:\Omega \rcirclearrowleft$ is the left shift map.
 Using the topological entropy, the previous results can be summarized as
\begin{equation}\label{eq:equality1}
\fbox{$h_{\top}(\sigma,K_p)=h_{\top}(\sigma,\Omega)-I(p), \; \forall p\in\R$},
\end{equation}
where $K_{p}:=\{\omega\in\Omega:\lim_n S_n/n=p\}$ 
{
and $S_n(\omega):=\sum_{i=0}^{n-1}\omega_i,\; \omega\in\Omega$.}
%
%
%
%
%
\subsection{Multifractal analysis of Birkhoff averages of H\"older continuous functions}
Most of the first results in the multifractal analysis of dynamical systems can be reduced to the results in the multifractal analysis of Birkhoff averages of certain potentials.
For simplicity, suppose $\Omega=\A^{\Z_+}$ is the set of all infinite
sequences $\omega=(\omega_0,\omega_1,\ldots)$ in a finite alphabet $\A$, equipped with a metric $\rho$, generating the product topology. 
Denote by $\sigma:\Omega\to\Omega$ the left shift on $\Omega$. 
The \mybf{topological pressure} of a continuous function  $\phi:\Omega\to\R$ is defined as
$$
P(\phi)=\lim_{n\to\infty}\frac 1n\log
\sum_{(a_0,\ldots,a_{n-1})\in A^n}
\Bigl( \sup_{\omega\in [a_0^{n-1}]}\exp(S_n\phi(\omega))\Bigr),
\quad
S_n\varphi(\omega) =\sum_{k=0}^{n-1}\phi(\sigma^k\omega),
$$
where $[a_0^{n-1}]:=\{\omega\in\Omega: \omega_i=a_i, \;i=0,...,n-1\}$ is the cylinder of the length $n$.
It is well-known that the pressure function $q\mapsto P_q(\phi):=P(q\phi)$, is convex. Moreover, by the celebrated result of D.Ruelle the pressure function for a \mybf{H\"older continuous} potential $\phi$
is \mybf{real-analytic}.

%
%

It turns out that the real-analyticity of the pressure function immediately gives us the
 Large Deviations Principle for ergodic averages $X_n(\omega)=S_n\phi(\omega)$.
More specifically, one has the following result\footnote{This is not the most general form of the Gartner-Ellis theorem.}
\begin{theorem}[Gartner-Ellis]\label{Gartner-Ellis}
 Suppose $\{X_n\}$ is a sequence of real-valued random variables defined
on a common probability space
$(\Omega,\mathcal F,m)$. Assume that for every $q\in\R$,  the logarithmic moment generating function, defined as the limit
$$
\Lambda(q) \triangleq \lim _{n \rightarrow \infty} \frac{1}{n} \log\mathbb E_m e^{ qX_n}=
\lim _{n \rightarrow \infty} \frac{1}{n} \log\int e^{ qX_n(\omega)} m(d\omega)
$$
exists and is finite. Denote by $\Lambda^*$ the Legendre transform (convex dual) of $\Lambda$, defined by
$$
\Lambda^*(\alpha) = \sup_{q\in \mathbb R}\Bigl(\alpha q- \Lambda(q)\Bigr).
$$
Then
\begin{itemize}
\item[(a)]  for any closed set $F\subset \R$,
$$
\limsup _{n \rightarrow \infty} \frac{1}{n} \log m\left(\left\{\omega\in\Omega:\  \frac 1nX_n(\omega)\in F\right\}\right)
 \leq-\inf _{\alpha \in F} \Lambda^*(\alpha) ;
$$

\item[(b)] if, furthermore $\Lambda(q)$ is differentiable on $\R$, for any open set $G\subset\R$,
$$\liminf _{n \rightarrow \infty} \frac{1}{n} \log
m\left(\left\{\omega\in\Omega:\  \frac 1nX_n(\omega)\in G\right\}\right)
 \geq-\inf _{\alpha \in G} \Lambda^*(\alpha).
 $$
\end{itemize}
\end{theorem}
  In order to apply the Gartner-Ellis theorem to ergodic averages $X_n(\omega)=S_n\phi(\omega)$
  it is necessary to introduce the reference measure $m$.
A good choice would be the {\it uniform } or  the {\it measure of maximal entropy} for $\Omega$. Namely, let
 $l=|\A|$ the number of different letters in alphabet $\A$, and let $m=\rho^{\Z+}_N$ be the product of uniform measures $\rho_N$ on $\A$. Then one immediately concludes that
 $$
 \Lambda_\phi(q) = P_\phi(q)-\log l= P(q\phi)-\log l.
 $$
 Therefore, for H\"older continuous functions $\phi$'s, the logarithmic moment generating function
 $\Lambda_\phi(q)$ is also
 real-analytic, and hence the Gartner-Ellis theorem is applicable. The so-called large deviations
 rate function is then
 $$\I_\phi(\alpha)=\Lambda^*_\phi(\alpha)=\sup_{q\in\R}\Bigl[q\alpha -\Lambda_\phi(q)\Bigr]
 =\sup_{q\in\R}\Bigl[q\alpha -P_\phi(q)+\log N\Bigr]=\log l+
P^*_\phi(\alpha).
 $$ Note also that $h_{\top}(\sigma,\Omega)=\log l$. Hence,
 $$
  \I_\phi(\alpha) =h_\top(\sigma,\Omega)+P_\phi^*(\alpha)
 $$

 \begin{theorem}[\cite{Barreira-Pesin-Schmeling1997}]\label{mfa:smooth} 
 Suppose $\Omega=A^{\Z_+}$ and $\phi:\Omega\to\R$ is H\"older continuous.
 For $\alpha\in\R$, consider the set
 $$
 K_\alpha=\left\{\omega\in\Omega:\
 \lim_{n\to\infty} \frac 1n \sum_{k=0}^{n-1}\phi(\sigma^k\omega)=\alpha\right\}.
 $$
 Then there exist $\underline\alpha,\overline\alpha\in\R$ such that
 \begin{itemize}
 \item for every $\alpha\in (\underline\alpha,\overline\alpha)$, $K_\alpha\ne\varnothing$  and
 $$\dim_{H}(K_\alpha) =\frac 1{\log l}\inf_{q}\Bigl( P(q\phi)-q\alpha\Bigr)
 =-\frac 1{\log l}\sup_{q}\Bigl( q\alpha -P(q\phi)\Bigr)=-\frac{P^*_\phi(\alpha)}{\log l}.
 $$
 \item the multifractal spectrum
 $$
 \alpha\mapsto\mathcal D_\phi(\alpha):= \dim_{H}(K_\alpha)
 $$
  is a real-analytic function of $\alpha$ on $( \underline\alpha,\overline\alpha)$.
 \end{itemize}
 \end{theorem}

Again, if we consider the topological entropy $h_{\top}(\sigma,\cdot)$ instead of the Hausdorff dimension
$\dim_{H}(\cdot)$, taking into account, that $h_{\top}(\sigma,\Omega)=\log l$, we obtain that the large deviations rate function and the topological entropy of level sets of ergodic averages are related
$$
h_{\top}(\sigma,K_\alpha)=h_{\top}(\sigma,\Omega)-\I_\phi(\alpha).
$$



It was observed in \cite{Takens-Verbitskiy1999}, that the result of Theorem \ref{mfa:smooth} can be extended to a much larger class of dynamical systems and observables $\phi$ with the key property being that for every $q\in\R$, there is a unique equilibrium
state $\mu_q$  for $q\phi$. To ensure uniqueness, one has to make two types of assumptions:
on the map $f:\Omega\to\Omega$ and the potential (observable) $\phi:\Omega\to \R$.
For example, in \cite{Takens-Verbitskiy1999} it was assumed that
\begin{itemize}
\item the map $f:\Omega\to\Omega$
is an \emph{expansive } homeomorphism with  \emph{specification} property
\item
the continuous  observable $\phi:\Omega\to\R$
is in the \emph{Bowen class}, i.e., for all $\epsilon>0$,
$$
\sup_{d_n(\omega,\omega')<\epsilon}\Bigl| (S_n\phi)(\omega)-(S_n\phi)(\omega')\Bigr|\le K(\epsilon)<\infty
$$
where the supremum is taken over all $\omega,\omega'$ such that
$$
d_n(\omega,\omega')=\max_{k=0,\ldots,n-1} d( f^k(\omega),f^k(\omega'))<\epsilon.
$$
\end{itemize}
Under these conditions, the Gartner-Ellis theorem holds for ergodic sums
$X_n(\omega)=S_n\phi(\omega)$ with a
differentiable  logarithmic moment generating function
$\Lambda(q)$, and the multifractal spectrum for ergodic averages is given by the Legendre transform
of the pressure function, i.e., for all $\alpha\in (\underline\alpha,\overline\alpha)$,
$$
h_\top(f,K_\alpha) =\inf_{q}\Bigl( P(q\phi)-q\alpha\Bigr)
 =-{P^*_\phi(\alpha)}. 
 $$ Furthermore, the logarithmic moment-generating function $\Lambda$ is related to the pressure function $P_{\phi}$ by $\Lambda(q)=P_{\phi}(q)-h_{\top}(f,\Omega),\; q\in\R$. 

%
%


\subsection{Multifractal analysis of Birkhoff averages of continuous functions}
If one considers simply continuous observables $\phi:\Omega\to\R$, it is no longer true
that  there is a unique equilibrium state $\mu_q$ for the potential $q\phi$ for all $q\in\R$.
Therefore, a different approach is required.

In the case of symbolic systems ($\Omega=A^{\Z_+}$, $\sigma:\Omega\to\Omega$ is  the left-shift),
the first results were obtained independently by Fan and Feng \cites{Fan1,Fan1996} and Olivier \cites{Olivier1, Olivier2}. The methods are quite different: Olivier relied on the density of H\"older continuous function
in the space of all continuous functions. Fan and Feng constructed the so-called Moran fractal.
The approach of Fan and Feng turned out to be very suitable for generalization to abstract compact spaces.
In \cite{Takens-Verbitskiy2003} the following variational principle for multifractal spectra has been obtained:

 \begin{theorem}[\cite{Takens-Verbitskiy2003}]\label{mfa:general} Suppose $T:\Omega\to\Omega$ is
 a continuous transformation of a compact metric space $(\Omega,d)$ with the specification property
 and $\phi:\Omega\to\R$ is a continuous observable.
 Then there exist $\underline\alpha,\overline\alpha\in\R$ such that
for every $\alpha\in (\underline\alpha,\overline\alpha)$, $K_\alpha\ne\varnothing$,  and
 $$h_\top(T,K_\alpha) =\sup_{\mu\in\mathcal M_T(\Omega)}\Bigl\{ h_\mu(T):\quad \int \phi d\mu=\alpha\Bigr\},
 $$
 where supremum is taken over $\mathcal M(\Omega,T)$ of all $T$-invariant measures on $\Omega$, and $h_\mu(T)$ is the
 Kolmogorov--Sinai entropy of $\mu$.
 \end{theorem}

%

For such systems, the Large Deviation Principle has been established in 1990 by Young \cite{Young}.
We will not recall the most general form of the result, but reformulate it in a such form
that the relation to the
multifractal formalism becomes apparent  immediately.

 \begin{theorem}[\cite {Young}, Theorem 1]\label{ldp:general}
Suppose  $m$ is a not necessarily invariant  Ahlfors-Bowen  reference measure on $\Omega$.
Suppose $\varphi:\Omega\to \R$ is continuous. Then for all $c\in \R$, one has
$$
\limsup_{n \rightarrow \infty} \frac{1}{n} \log m\left\{\frac{1}{n} S_n \varphi \ge c\right\}
\le \sup_{\mu\in\mathcal M_T(\Omega)} \left\{h_\mu(T):   \int \varphi d \mu \geq c\right\} -h_{\top}(T,\Omega).
$$
If, furthermore, $T$ has a specification property, then for all $c\in\R$, one has
$$
\liminf_{n \rightarrow \infty} \frac{1}{n} \log m\left\{\frac{1}{n} S_n \varphi >c \right\}
\ge
\sup_{\mu\in\mathcal M_T(\Omega)} \left\{h_\mu(T):   \int \varphi d \mu > c\right\} -h_{\top}(T,\Omega).
$$
\end{theorem}
This result allows us to say, that under the above conditions, the sequence of random variables $\{X_n(\omega)=(S_n\phi)(\omega)\}$ satisfies the Large Deviation Principle with the rate function
$$ I_\phi(\alpha) = h_{\top}(T,\Omega) - \sup_{\mu\in\mathcal M_T(\Omega)}\Bigl\{ h_\mu(T):\quad \int \phi d\mu=\alpha\Bigr\},
$$
As in the cases considered above, combining the results of Theorems \ref{mfa:general}  and  \ref{ldp:general}, one concludes that for topological dynamical systems with specification the large deviations rate
function $I_{\phi}$ and the multifractal spectrum $\e_{\phi}(\alpha):=h_{\top}(T,K_{\alpha})$ are related by
$$
\mathcal E_\phi(\alpha) = h_\top(T,\Omega)-I_\phi(\alpha).
$$
%
In fact, the methods  developed by Young are instrumental in the proof of Theorem \ref{mfa:general}  in \cite{Takens-Verbitskiy2003}.

\subsection{Multifractal analysis beyond Birkhoff sums}
A number of multifractal results have been obtained for the level sets of the form
$$
K_\alpha=\Bigl\{\omega\in\Omega:\, \lim_{n\to\infty} \frac 1nX_n(\omega)=\alpha\Bigr\},
$$
where $X_n$'s are not necessarily Birkhoff sums of the form $S_n\phi$, for some $\phi$.
For example, such \emph{non-additive} observable arise naturally in the multifractal analysis of pointwise Lyapunov exponents for the products of matrices. Namely, suppose a map $M:\Omega\to M_d(\R)$ is given where $M_d(\R)$ is the algebra of $d\times d$ real matrices. We can associate to every $x\in\Omega$ and $n\in\N$, a random variable $X_{n}(x):=\log||M(T^{n-1}x)M(T^{n-2}x)\cdots M(x)||$ provided that $||M(T^{n-1}x)M(T^{n-2}x)\cdots M(x)||\neq 0$, where $||\cdot||$ is a norm on $M_d(\R)$. Clearly, $\{X_n\}_{n\in\N}$ is not an additive sequence if $d\geq 2$.  

A non-additive sequence of potentials can also be encountered if one considers weighted ergodic sums $\{S_n^{w}f\}$ of a continuous potential $f$ on the space $\Omega$, i.e., $S_n^{w}f:=\sum_{i=0}^{n-1}w_if\circ T^i,\; n\in\Z_+$ and $w_n\in\R$ are weights \cite{Fan2021}. It is easy to see that $\{S_n^{w}f\}$ is not additive unless all the weights $w_n,\; n\in\Z_+$ are the same. 

Another example of non-additive sequences arises in studying pointwise entropies of probability measures. Suppose $\mu$ is a fully supported probability measure on $\Omega$. For a fixed $\delta>0$, define a sequence of random variables as $X_n(x):=-\log\mu (B_n(x,\delta)),\;n\in\N$, $x\in\Omega$, where $B_n(x,\delta):=\{y\in\Omega:d(T^iy, T^ix)<\delta,\;\; i=0,\dots, n-1\}$, $n\geq 1$. 

We shall discuss these and other examples in greater detail in Section \ref{examples} and show how our results can be applied in these cases.




\section{Preliminaries and Notations}

In this section, we assume that we are given a compact metric space $(\Omega,d)$ equipped with a continuous transformation $T:\Omega\to\Omega$.
We denote by $M(\Omega), M(\Omega, T)$ the sets of all Borel probability measures and $T$-invariant Borel probability measures on $\Omega$, respectively.
For $x,y\in\Omega$ and $n\in\N$, the Bowen metric $d_n$ is defined as $d_n(x,y):=\max_{0\leq i\leq n-1} d(T^ix,T^iy)$.
By $B_n(x,\varepsilon)$ we denote an open ball of radius $\varepsilon> 0$ in the metric $d_n$ centred at $x\in \Omega$, i.e., $B_n(x,\varepsilon):=\{y\in\Omega: d_n(x,y)<\varepsilon\}$.

\subsection{ Topological Entropy}\

The notion of topological entropy of non-compact and non-invariant sets has been introduced by Bowen \cite{Bowen} in 1973.
This paper uses an equivalent definition of Bowen's topological entropy given in \cite{Pesin}.
For $Z\subset \Omega$, and all $t\in \R,\; \epsilon>0$,  and  $N\in\N$, let
\begin{equation}\label{non-centered def of top entropy}
   m(Z,t,\epsilon,N):=\inf \left \{\sum_{i=1}^{\infty} e^{-n_it}: Z\subset \bigcup_{i=1}^{\infty} B_{n_i}(x_i,\epsilon),\; n_i\geq N \right \} .
\end{equation}
By the standard convention, $m(\emptyset,t,\epsilon,N)=0$, for all $ t,\epsilon, N$.
Obviously, $m(Z,t,\epsilon,N)$ does not decrease on $N$, thus we can define
$$m(Z,t,\epsilon):=\lim_{N\to\infty}m(Z,t,\epsilon,N).$$
One can show \cite{Pesin} that $m(\cdot,t,\epsilon)$ is an outer measure with properties similar to those of $t$-dimensional Hausdorff outer measure.
In particular, there exists a critical value $t'\in\R$ such that
 \begin{equation*}
 m(Z,t',\epsilon)=\begin{cases}  +\infty, \quad ~~~ \text{if}~  t'<t, \\
  0, \quad   ~~~~~\text{if}~  t'>t.
  \end{cases}
 \end{equation*}
 We denote this critical value by $h_{\top}(T,Z,\epsilon)$.
 Thus $h_{\top}(T,Z,\epsilon)=\inf\{t\in\R: m(Z,t,\epsilon)=0\} =\sup \{ t\in\R: m(Z,t,\epsilon)=+\infty\} $.
 Furthermore, since $h_{\top}(T,Z,\epsilon)$ is monotonic in $\epsilon$,
 we can define the \mybf{topological entropy} of $Z$ by 
$$h_{\top}(T,Z):=\lim\limits_{\epsilon\to 0+}h_{\top}(T,Z,\epsilon).$$
It should be stressed that the set $Z$ is not assumed to be compact nor $T$-invariant.
Finally, by our convention, $h_{\top}(T,\emptyset)=-\infty$, and $h_{\top}(T,Z)\geq 0$ for all non-empty subsets $Z\subset \Omega$.
We summarize the basic properties of the topological entropy in the following remark.
\begin{remark}[\cite{Pesin}] \label{Pesin's-book}
 \begin{itemize}
        \item [(i) ] \textit {Monotonicity: } if $Z_1\subset Z_2$, then $h_{\top}(T,Z_1)\leq h_{\top}(T,Z_2)$;
        \item[ (ii)] \textit {Countable stability: } if $Z=\bigcup_n Z_n$, then $h_{\top}(T,Z)=\sup_{n} h_{\top}(T,Z_n)$.
    \end{itemize}
\end{remark}

\begin{remark}\label{Hausdorff dim and top entropy}
The notions of topological entropy and Hausdorff dimension coincide if the underlying space $\Omega$ is a symbolic space. Suppose $\Omega =\{1,2, \ldots, l\}^{\mathbb{Z}_+}$ with integer $l \geq 2$, and the metric on $\Omega$ is defined as
 \[
 d(x, y)=\frac{1}{l^k},  \   \ \text{where} \  \  k=\min \{i\geq 0: x_i \neq y_i\},
 \]
then $Z \subseteq \Omega$, $\dim_H (Z) = \frac{h_{\top} (T, Z)}{\log l}$.
\end{remark}

\subsection{Local (pointwise) entropies}\

Consider a Borel probability measure $\mu$ on the metric space $(\Omega,d)$, and define the lower and upper pointwise entropies of $\mu$ at $x \in \Omega$ as
\begin{equation*}\label{LUE}
  \underline{h}_{\mu} (T, x) : =\lim_{\varepsilon \to 0} \liminf_{n \to \infty} -\frac{1}{n} \log \mu(B_n(x, \varepsilon)),\;\;\; \overline{h}_{\mu} (T, x) : =\lim_{\varepsilon \to 0} \limsup_{n \to \infty} -\frac{1}{n} \log \mu(B_n(x, \varepsilon)).
\end{equation*}

Note that the limit in $\varepsilon$ exists due to the monotonicity. 
D.J. Feng and W. Huang \cite{Feng-Huang2012} extended the notion of measure-theoretic entropy to non-invariant measures by motivating the result of M. Brin and A. Katok. In fact, for a probability measure $\mu\in \M(\Omega)$, D.J. Feng and W. Huang \cite{Feng-Huang2012} defined the \textit{upper} and \textit{lower measure-theoretic entropies} of $\mu$ relative to the transformation $T$ by 
$$\underline{h}_{\mu}(T):=\int_{\Omega}\underline{h}_{\mu}(T,x)\mu(dx) \text{ and } \overline{h}_{\mu}(T):=\int_{\Omega}\overline{h}_{\mu}(T,x)\mu(dx).
$$
Note that the above formulas are consistent with the classical notion of measure-theoretic entropy due to the Brin-Katok theorem \cite{Brin}. In fact, if $\mu\in \M(\Omega,T)$, then $\underline{h}_{\mu}(T)=\overline{h}_{\mu}(T)=h_{\mu}(T)$, where $h_{\mu}(T)$ is the measure-theoretic entropy of $\mu$ with respect to $T$. In \cite{Feng-Huang2012}, the authors also extended the classical variational principle to the non-invariant setting. More precisely, they showed for a compact, but not necessarily invariant set $K\subseteq\Omega$ that  
\begin{equation}\label{VP for topological entropies of non-invariant sets}
    h_{\top}(T,K)=\sup\{\underline{h}_{\mu}(T):\mu\in\M(\Omega),\; \mu(K)=1\}.
\end{equation}




The next theorem is a dynamical analogue of the Mass distribution principle and is very useful for the estimation of topological entropies of various sets.  
\begin{theorem}[\cite{Ma-Wen}]\label{entropy distribution principle}
Let $\mu$ be a Borel probability measure on $\Omega$, and $E$ be a Borel subset of $\Omega$.
\begin{itemize}
    \item[(1)]
     If~ $\underline{h}_{\mu}(T, x)\le t$~ for ~ some ~ $t<\infty$ ~ and ~all ~$x\in E$,  then $h_{\top}(T, E)\le t$.

    \item[(2)]
    If ~ $\mu(E)>0,$~ and ~ $\underline{h}_{\mu}(T, x)\ge t$~ for ~ some ~ $t>0$ ~ and ~ all ~$x\in E$, then $h_{\top}(T, E) \ge t$.
\end{itemize}
\end{theorem}

%
\subsection{Ahlfors-Bowen measures}

\begin{defn}
We say  that a not necessarily invariant Borel probability measure $m$ on $\Omega$ is \mybf{Ahlfors-Bowen}
if there exists $h>0$ such that for every $\epsilon>0$ there exists a finite positive constant, $C(\eps)$ such that for every $n\ge 1$ and all $\omega\in\Omega$, one has
\begin{equation}\label{Ahlfors-Bowen-measure}
\frac 1{C(\eps)} e^{-nh}\le m\bigr( B_n(\omega,\epsilon)\bigr)\le
{C(\eps)} e^{-nh}.
\end{equation}
\end{defn}
It is easy to see from the definition that an Ahlfors-Bowen measure has a very simple spectrum of local entropies. More specifically, $h_{m}(T,\omega)=h$ for all $\omega\in\Omega$. Furthermore, one has the following simple lemma for the constant $h$. 


\begin{lemma}\label{A-B_constant_and_top_entropy}
Suppose $m$ is Ahlfors-Bowen measure for $T$ on $\Omega$, then $h=h_\top(T,\Omega)$.
\end{lemma}
\begin{proof} 
Since $m$ is an Ahlfors-Bowen measure, for any $\omega \in \Omega$, one has $h_{m} (T, \omega) = h$.
Therefore, by Theorem \ref{entropy distribution principle}, we immediately conclude that $h_{\top}(T, \Omega) = h$.
\end{proof}


\subsection{Conditions on the underlying transformation}\label{expansiveness and specification}\

In applications, topological dynamical systems (TDS) usually have some sort of expansivity and mixing properties. Therefore, there is a vast amount of literature devoted to studying dynamical systems with these properties. We also impose expansivity and mixing properties on our underlying transformation to prove some of our results.   

\begin{defn}
  A continuous transformation $T: \Omega \to \Omega$ is called \mybf{expansive} if there exists a constant $\rho >0$ such that if
  $$
  d(T^n (x), T^n (y) ) < \rho\  \text{for all non-negative integer } n, \text{then}\ \;  x=y.
  $$
  The maximal $\rho$ with such a property is called the expansive constant.
\end{defn}

A slightly stronger version of topological mixing is topological exactness. Topological exactness of the transformation $T$ is that for every $\epsilon>0$ and $x\in\Omega$ there is $N:=N(x,\epsilon)\in \N$  such that $T^N(B(x,\epsilon))=\Omega$. Basically, by using the compactness of the underlying space $\Omega$, one can uniformise the constant $N$ with $x\in\Omega$.  A relatively stronger version of topological exactness which we call strong topological exactness also gives uniform control over the Bowen balls of the transformation. 
\begin{defn}
We call a continuous transformation $T:\Omega\to\Omega$ \mybf{strongly topologically exact} if for any $\epsilon>0$ there is a natural number $M_1\in \N$ such that for all $n\in\N$ and $x\in\Omega$, 
\begin{equation}\label{strong topological exactness}
T^{n+M_1}(B_n(x,\epsilon))=\Omega.
\end{equation}
\end{defn}
Note that topologically mixing subshifts of finite type are strongly topologically exact.
P.Walters (\cite{Walters2000}) showed that an expansive strongly topologically exact transformation has a weak form of the specification property, namely, one has the following theorem:
%
\begin{theorem}[\cite{Walters2000}]\label{Walters theorem}
Let $T: \Omega \to  \Omega$ be an expansive strongly topologically exact transformation, then  $T$ satisfies the \textit{weak specification condition}, i.e., for any $\varepsilon>0$ there exists a natural number $M_2\in\N$ such that for for all  $x, x' \in \Omega$ and every $n_1, n_2 \in\N$, there exists $w \in \Omega$
with~$d(T^i w, T^i x) <\varepsilon, 0 \leq i \leq n_1-1$~and $d(T^{n_1+M_2+j} w, T^j x') <\varepsilon,\; 0 \leq j \leq n_2-1$. Note that the latter condition is equivalent to the following: for every $\epsilon>0$ there exists an integer $M_2>0$ such that for all $x, x' \in \Omega$ and $n_1, n_2 \in\N$, one has
\begin{equation}\label{weak specification condition}
B_{n_1}(x,\epsilon)\cap T^{-(n_1+M_2)}(B_{n_2}(x',\epsilon))\neq \emptyset.
\end{equation}

\begin{remark}
If $T: \Omega \to \Omega$ is an expansive homeomorphism with weak specification property, then the measure of maximal entropy for $T$ is Ahlfors-Bowen \cite[Theorem 4.6]{Walters2000}. 
\end{remark}
\end{theorem}

\subsection{Class of random variables}\

We say  that the sequence $\{X_n\}_{n\geq 1}$  of real-valued random variables (observables) on $\Omega$ is  \mybf{subadditive} if for every $m,n \geq 1$,
$$X_{n+m} \leq X_n +X_m \circ T^n .$$
We call a sequence $\{X_n\}_{n\geq 1}$  of random variables on $\Omega$ \mybf{weakly almost additive} (or simply, almost additive),
if there are constants $A_n \geq 0,\; n\in\N$ with $\lim_{n\to\infty}\frac{A_n}{n}=0$ such that for all $x \in \Omega$ and $n, m\in \mathbb{N}$,
\begin{equation}\label{almost additivity}
|X_{n+m}(x)- X_n(x)- X_m(T^n x)  | \leq A_n .
\end{equation}
Recently, Cuneo \cite[Theorem 1.2]{Cuneo2020} showed that weakly almost additive sequences of continuous functions can be approximated 
uniformly by Birkhoff's sums of continuous functions.


For a function $H:\Omega\to\R$ and $n\in\N$, $\delta>0$, we denote the $(n,\delta)$-variation of the function $H$ by $v_{n,\delta}(H):=\sup\limits_{d_n (y,z) \leq \delta} | H (y)- H (z)| $.
\begin{defn}
A sequence $\{X_n\}$ of random variables on $\Omega$ satisfies the \mybf{Bowen condition},
if there exists $\delta>0$ such that $\sup\limits_{n\in\N} v_{n,\delta}(X_n)<+\infty$. We say that the sequence $\{X_n\}$ satisfies the \mybf{weak Bowen condition} if there exists $\delta>0$ such that
\begin{equation}\label{weak Bowen condition}
\lim_{n\to\infty} \frac{v_{n,\delta}(X_n)}{n}=0.
\end{equation}
\end{defn}
It should be mentioned that the continuity of the random variables $X_n, \; n\geq 1$, is not assumed in the definition of the weak Bowen condition, and neither the Bowen condition nor the weak Bowen condition implies the continuity of each individual random variable $X_n$. However, if a sequence of functions $\{X_n\}$ satisfies the weak Bowen condition, then each function $X_n$ is bounded.
Indeed, since $\Omega$ is compact, for every $n$, there exists a finite cover $\{B_n({\omega_i,\delta/2})\}_{i}$,
and the observable $X_n$ is a bounded function on each ball $B_n({\omega_i,\delta/2})$ in the cover. 

\begin{remark}\label{expansivity important property}
One can observe that if $T$ is \textit{expansive}, then the sequence of ergodic sums of any continuous potential $\psi:\Omega\to\R$ satisfies weak Bowen's condition.
Namely, if $X_n=S_n\psi=\sum_{i=0}^{n-1}\psi\circ T^i$, then $v_{n,\delta}(X_n)\leq \sum_{i=1}^{n}v_{i,\delta}(\psi)$, and thus 
\[
\lim_{n\to\infty}\frac{v_{n,\delta}(X_n)}{n}\leq \lim_{n\to\infty}\frac{\sum_{i=1}^{n}v_{i,\delta}(\psi)}{n}=0
\]
by Stolz-Cesaro's theorem, since $\psi:\Omega\to\R$ is continuous and hence $\lim\limits_{n\to\infty} v_{n,\delta}(\psi)=0$, where positive number $\delta$ is smaller than the expansive constant of $T$.
\end{remark}

\section{Large Deviations}
\subsection{Some concepts from convex analysis} In this section, we consider functions with values in $\overline{\R}:=\R\cup\{\pm \infty\}$. A convex function $\phi:\R\to\overline{\R}$ is called a \textbf{proper convex function} if $\phi(t)>-\infty$ for all $t\in\R$ and $\phi\nequiv +\infty$. The \mybf{essential domain} of the convex function $\phi$ is $\dom(\phi):=\{t\in\R:\phi(t)<+\infty\}$. Note that  $\dom(\phi)\subset \R$ is a convex set.

Suppose  $\phi:\R\to (-\infty,+\infty]$ is a proper convex function and assume that $\phi$ is finite at some $t\in\R$. A real number $t^*$ is a \textbf{subgradient} of $\phi$ at $t$ if for all $s\in\R$, $\phi(t)+t^*(s-t)\leq \phi(s)$ holds.  Denote the set of all subgradients of $\phi$ at $t$ by  $\partial \phi(t)$. It is easy to check that  $\partial\phi(t)=\R\bigcap \Big[\sup_{s>0}\frac{\phi(t)-\phi(t-s)}{s},\; \inf_{s>0}\frac{\phi(t+s)-\phi(t)}{s}  \Big]=\R\bigcap[\phi'_{-}(t),\phi'_+(t)]$, where $\phi'_{-}(t),\; \phi'_{+}(t)\in \overline{\R}$ are left and right derivatives of $\phi$ at $t$. We define the domain of the multivalued mapping $\partial\phi: t\in \dom(\phi)\mapsto \partial \phi(t)\subset\R$ by $\dom(\partial \phi):=\{t\in \dom(\phi): \partial \varphi (t)\neq\emptyset\}$. Note that in general,  $\dom(\partial \phi)$ is not a convex set. 
\begin{defn}
A proper convex function $\phi:\R\to\overline{\R}$ is called \mybf{essentially strictly convex} if $\phi$ is strictly convex on every convex subset of $\dom(\partial \phi)$. Note that a convex function $\phi$ is called \textit{strictly convex} in a convex set $D\subset \R$, if for every $t_1,t_2\in D$ with $t_1\neq t_2$ and $\lambda\in(0,1)$, $\phi(\lambda t_1+(1-\lambda)t_2)<\lambda\phi(t_1)+(1-\lambda)\phi(t_2)$. 
\end{defn}
It is a well-known fact that finite convex functions are very close to differentiable functions. However, if it comes to extended convex functions, then one needs to be more careful about the essential domain of the function.  
\begin{defn}
A proper convex function $\phi:\R\to\overline{\R}$ is \mybf{essentially smooth} if it satisfies the following conditions:
\begin{itemize}
\item[(a)] $\dom (\phi)^{\circ}\neq \emptyset$, where $A^{\circ}$ denotes the interior of a set $A$;
\item[(b)] $\phi$ is differentiable on $\dom (\phi)^{\circ}$;
\item[(c)] $\phi$ is steep, namely, $\lim _{n \rightarrow \infty}\left|\phi\left(t_n\right)\right|=\infty$ whenever $\left\{t_n\right\}$ is a sequence in $\dom (\phi)^{\circ}$ converging to a boundary point of $\dom (\phi)^{\circ}$.
\end{itemize}
\end{defn}
An extended, not necessarily convex, function $\phi:\R\to\overline{\R}$ is \mybf{lower semicontinuous} if for all $\lambda\in \R$,  $\{t\in\R:\phi(t)\leq \lambda\}$ is a closed subset of $\R$. Equivalently, a function $\phi$ is lower semicontinuous if and only if $\liminf\limits_{x\to x_0}\phi(x)\geq \phi(x_0)$ for all $x_0\in\R$. Note that a lower semicontinuous function achieves its minimum in any compact set.

The following theorem shows that essentially strictly convex and essentially smooth convex functions are intimately related in terms of the Legendre transform. 

\begin{theorem}[\cite{Rockafellar-Book}]\label{the theorem in Rockafellar book}
A lower semicontinuous proper convex function $\phi$ is essentially strictly convex if and only if its Legendre transform, $\phi^*(t):=\sup\limits_{q\in \R} \{tq-\phi(q)\},\; t\in\R$, is essentially smooth.
\end{theorem}

 \subsection{Large deviations}
\begin{defn}\label{ldp rate function definition}
An extended function $I:\R\to[0,+\infty]$ is called a \mybf{ large deviations rate function} (or simply a rate function) if it is lower semicontinuous.
A rate function $I$ is called \textbf{good} if the sub-level sets, $\{t:I(t)\leq \lambda\}\subset\R$ are compact for all $\lambda\in\R$. For $E\subset\R$, $I(E)$ abbreviates $\inf\limits_{t\in E} I(t)$ i.e., $I(E):=\inf\limits_{t\in E} I(t)$. 
\end{defn}
Note that a good rate function achieves its minimum in any closed set. 

Suppose $(\Omega, \mathcal{F}, \mu)$ is a probability space
and $X_n:\Omega\to\R$, $n\ge 1$, are random variables. 

\begin{defn} The sequence of random variables $\{\frac{1}{n}X_n\}_{n=1}^{+\infty}$
satisfies the \mybf{large deviations principle} (LDP) with a \mybf{rate function $I$} if
for all Borel sets $E\subset \R$,
$$- I(E^{\circ}) \leq \liminf _{n \rightarrow \infty} \frac
1n\log \mu\Big(\Big\{\omega: \ \frac{X_n(\omega)}{n}\in E\Big\}\Big)
 \leq
 \limsup _{n \rightarrow \infty} \frac
1n\log \mu\Big(\Big\{\omega: \ \frac{X_n(\omega)}{n}\in E\Big\}\Big) \leq-I(\overline{E}),
$$
where $\overline{E}$, $E^{\circ}$ denote the closure and the interior of $E$,
respectively.
\end{defn}
One of the most convenient ways to  establish the validity of the LDP is provided by the Gartner-Ellis theorem which we have stated in a partial case in Theorem \ref{Gartner-Ellis}. 
 Denote the normalised log-moment generating function of $X_n$ by $\phi_n$, i.e., for $t\in\R$ 
$$\phi_n(t):=\frac{1}{n}\log \int\limits_{\Omega}e^{tX_n}d\mu.$$ 
Note that $\phi_n$ may take value $+\infty$. Define the log-moment generating function of the sequence $X:=\{X_n\}$ by $\phi_X(t):=\limsup\limits_{n\to\infty}\phi_n(t),\; t\in\R$. Then the Gartner-Ellis theorem states that if for all $t\in\R$, the limit, $\lim\limits_{n\to\infty}\phi_n(t)$ exists and is finite, and the log-moment generating function $\phi_X$ is essentially smooth, then 
$\{\frac{1}{n}X_n\}_{n\in\N}$ satisfies the LDP with a good rate function $I_X:=\phi_X^*$. An inverse statement to this in some sense is the celebrated Varadhan's lemma which we reformulate slightly for our purposes below. 
One can find a general version of the theorem in \cite[Theorem 4.3.1]{Dembo-Zeitouni-Book}.

\begin{theorem} [\cite{Dembo-Zeitouni-Book}]\label{Varadhan theorem}
Suppose that $\{\frac{1}{n}X_n\}_{n=1}^{+\infty}$ satisfies the LDP with a good rate function $I: \R \rightarrow[0, \infty]$, and 
$$
\limsup _{n \rightarrow 0} \frac 1n\log \int_{\Omega}e^{t X_n }d\mu=\phi_X(t)<\infty 
$$
hold for all $t\in\R$.
Then the limit $\lim\limits _{n \rightarrow 0} \frac 1n\log \E_{\mu}\left[e^{t X_n }\right]$ exists and
$$
\phi_X(t)=\lim\limits _{n \rightarrow 0} \frac 1n\log \int_{\Omega}e^{t X_n }d\mu=\sup_{q \in \R}\{tq-I(q)\}=I^*(t),\;\;\; \text{for all }\; t\in\R.
$$
\end{theorem}

Sometimes, a weak version of LDP is more convenient for purposes, since it is easier to establish than the full LDP. 
\begin{defn}
Let $(\Omega,\mathcal{F},\mu)$ be a probability space and $X_n:\Omega\to\R,\;\; n\geq 1$ be random variables. The sequence $\{\frac{1}{n}X_n\}_{n\in\N}$ satisfies the \textbf{weak Large Deviation Principle} (weak LDP) with a rate function $I$ if 
\begin{itemize}
\item[1)] for all compact set $F\Subset\R$, 
\[
\limsup_{n\to\infty}\frac{1}{n}\log\mu\Big(\left\{ \frac{1}{n}X_n\in F \right\}\Big)\leq -I(F)
\]
\item[2)] for all open set $G\subset\R$, 
\[
\liminf_{n\to\infty}\frac{1}{n}\log\mu\Big(\left\{ \frac{1}{n}X_n\in G \right\}\Big)\geq -I(G)
\]
\end{itemize}
\end{defn}
Note that the only difference between LDP and weak LDP is in the upper large deviation bound.

As one can see in the Gartner-Ellis theorem, the LDP upper bound requires fewer conditions than the lower bound. 
However, one still needs to make rather mild assumptions (for example, the existence of the limit log-moment generating function as in Theorem \ref{Gartner-Ellis}) to obtain the upper bound. Nevertheless, one can establish the upper \textit{weak large deviation}  bound without any additional assumption about the sequence as the following folklore theorem states.
%
\begin{theorem}\label{Gartner-Ellis UB weak LDP}
Let $\{X_n\}$ be a sequence of random variables on a probability space $(\Omega,\mathcal{F},\mu)$. Then the sequence $\{\frac{1}{n}X_n\}_{n\in\N}$ satisfies the weak large deviations upper bound with the rate function $I_X:=\phi_X^*$, i.e., for any compact set $F\Subset \R$, the following holds
\begin{equation}
\limsup_{n\to\infty}\frac{1}{n}\log\mu\Big(\left\{ \frac{1}{n}X_n\in F \right\}\Big)\leq -I_X(F).
\end{equation}
\end{theorem}
   


%
%


\section{Main results}

The main goal of this paper is to obtain the Multifractal spectra for the sequence $\{\frac{1}{n}X_n\}$ under the assumption that the sequence satisfies LDP. As a byproduct, we also obtain a formula that relates the entropy spectra with the LDP rate function. Although we assume in our theorems that the reference measure satisfies the Ahlfors-Bowen condition, one should be able to generalise  the results to a multifractal measure $\tau$ with the smooth multifractal spectrum, i.e. $\alpha\mapsto h_{\top}(T,E_{\alpha})$ is a smooth function, where $E_{\alpha}:=\{\omega: h_{\tau}(\omega)=\alpha\}$.

 Assume that the ambient space $(\Omega,d)$ is compact, and we have a continuous transformation $T:\Omega\to\Omega$ with finite topological entropy, i.e., $h_{\top}(T,\Omega)<\infty$. It is also assumed that a sequence $X:=\{X_n\}$ of random variables on $\Omega$ and a reference measure $\nu\in\mathcal{M}(\Omega)$ satisfying Ahlfors-Bowen condition with the parameter $h=h_{\top}(T,\Omega)$ are given. 
It should be stressed that we are not assuming $\{X_n:n\in\N\}$ are continuous functions, and $\nu$ is $T$-invariant. We denote, for $\alpha\in\R$, the $\alpha$-level set of the sequence $\{\frac{1}{n} X_n\}$ by
$$ K_{\alpha}:=\{x\in\Omega: \lim\limits_{n\to\infty}\frac{1}{n} X_n(x)=\alpha\},
$$ and the domain of the multifractal spectrum by  $\La_X:=\{\alpha\in\R:K_{\alpha}\neq\emptyset\}$. We shall study the entropy spectrum $\mathcal{E}_{X}(\alpha):=h_{\top}(T,K_{\alpha})$. In the statements of our theorems, we shall assume some of the following main conditions:
\begin{itemize}
\item[(A1)] The sequence $\{X_n\}_n$ satisfies the weak Bowen condition;
\item[(A2)] The sequence $\{X_n\}_n$ is weakly almost additive;
\item[(A3)] $T: \Omega \to \Omega $ is an expansive strongly topologically exact transformation.
\end{itemize}
%
We now present our first main result that estimates the entropy spectrum $\mathcal{E}_{X}$ from above in terms of the LDP rate function.

\begin{theorema} Assume (A1), and the sequence $\Big\{\frac{1}{n}X_n\Big\}$ satisfies {\normalfont{the weak large deviations upper bound}} with a rate function $I_X:\R\to[0,+\infty]$, i.e.   for all compact $ F\Subset\R$ one has
    \begin{equation}\label{LDUB}
      \limsup_{n\to\infty}\frac{1}{n}\log\nu\Big\{\frac{1}{n}X_n\in F \Big\}\leq -I_X(F),
    \end{equation}
then for all $\alpha\in\La_X$, one has
    \begin{equation}\label{MFUB_from_LDUB}
    \e_X(\alpha)\leq h_{\top}(T,\Omega)-I_X(\alpha).
    \end{equation}
\end{theorema}    
%
%
%
Let $\phi_{n}(q):=\frac{1}{n}\log\E_{\nu}e^{qX_n}$ and $\phi_X:=\limsup\limits_{n\to\infty}\phi_{n}$  be the log and log-limit moment-generating functions of the random variable $X_n$, respectively. Notice that the moment-generating functions $\phi_n$ exist and are finite for all $n\in\N$ since the observables $X_n$ are bounded functions by the weak Bowen condition. The next result estimates the spectrum $\e_X$ from above in terms of the log-moment generating function $\phi_X$ of $X$.



%
\begin{theoremb}\label{result_1}
Assume (A1), then for all $\alpha\in\La_X$
\begin{equation}\label{ub with moment generating function}
\mathcal{E}_X(\alpha)\leq h_{\top}(T,\Omega)-\phi_X^*(\alpha),
\end{equation}
where $\phi_X^*$ is the Legendre-Fenchel transform of $\phi_X$.
\end{theoremb}
\begin{remark}
If one drops the condition (A1) in the above theorems, then the statements of the theorems will be false. Assume $\Omega:=\{0,1\}^{\Z_+}$, $T:\Omega\to\Omega$ is the left shift transformation, $x\in\Omega$ and $\kappa\in\R\setminus\{1\}$.  Set $X_n(\omega):=n$ if $\omega_i=x_i, \;i=0,1, \ldots, n^2-1 $, and $X_n(\omega):=\kappa n$ otherwise. Then $\phi_X(q)=\kappa q$ for every $q\in\R$, hence by Theorem \ref{Gartner-Ellis}, the sequence $\{X_n\}_{n\in\N}$ satisfies LDP with a good rate function $\phi_X^*$. However, one has $h_{\top}(T,K_1)=0>h_{\top}(T,\Omega)-\varphi_X^*(1)=\log2-\infty=-\infty$.
\end{remark}

\begin{remark}
We should note that one can prove the same statements as in Theorem A and Theorem B with the following condition which is slightly weaker than the weak Bowen condition: there exists $\delta>0$ such that
\[
\liminf\limits_{n\to\infty}\frac{v_{n,\delta}{(X_n)}}{n}=0.
\]
\end{remark}

One might expect equality in Theorem A and Theorem B, but it is easy to show by examples that the conditions of these theorems are ``too general" to prove such equality.  The subsequent theorem provides sufficient conditions for the equality to hold in Theorem A and Theorem B.

\begin{theoremc}\label{lb theorem}
Assume the conditions (A1), (A2) and (A3). If the sequence $\big\{\frac{1}{n} X_n\big\}_{n\in\N}$ satisfies LDP with an essentially strictly convex good rate function $I_X$,
then one has the following:
\begin{itemize}
\item[(i)] There exists extended real numbers $-\infty\leq\underline{\alpha}\leq\overline{\alpha}\leq+\infty$ such that $(\underline{\alpha},\overline{\alpha})\subset \L_X\subset [\underline{\alpha},\overline{\alpha}]$.
\item[(ii)] $\e_X(\alpha) =h_{\top}(T,\Omega) - I_X(\alpha)$ holds for all $\alpha\in (\underline{\alpha},\overline{\alpha})$.
\end{itemize}
\end{theoremc}

Finally, it is well-known that a result of multifractal formalism for random variables can be generalized straightforwardly to random vectors (see \cite{Takens-Verbitskiy2003}). In particular, one can directly extend the above theorems to sequences $\{\Vec{X}_n\}_{n\in N}$ of random vectors $\Vec{X}_n:\Omega\to \R^d$ with straightforward modification of the conditions (A1)-(A3). 

\section{Proofs of Theorem A and Theorem B: Upper bounds}\label{proof ub-section}
\subsection{A covering lemma} The Vitali covering lemma presented below plays a crucial role in establishing upper bounds and is also instrumental in proving a key lemma in the subsequent sections.
Therefore, we shall state and  prove it first for the sake of completeness.
\begin{lemma}[Vitali covering lemma]\label{Vitali covering lemma}
Let $(\Xi,d)$ be a metric space, and suppose $\mathfrak{F}$ is a non-empty finite collection of balls. Then there exists a disjoint subcollection $\mathfrak{F}'\subseteq \mathfrak{F}$, such that $\bigcup\limits_{B\in\mathfrak{F}} B\subset \bigcup\limits_{B'\in\mathfrak{F}'} 3B'$, where $3B'$ denotes a ball with the same centre as the ball $B'$, but with three times the radius of $B'$.
\end{lemma}
\begin{proof}
We shall iteratively construct an ascending chain $\{\mathfrak{F}_s\}$ of collections of pairwise disjoint balls contained in $\mathfrak{F}$ in the following way: assume that in $s$-th step ($s\in\N$),  we have constructed such a subcollection  $\mathfrak{F}_s$, ($\mathfrak{F}_1$ can contain only one ball, therefore, such a disjoint collection of the balls always exists) then we set $\mathfrak{F}_{s+1}=\mathfrak{F}_s$ if $B'\bigcap \bigcup\limits_{B\in\mathfrak{F}_s}B\neq \emptyset$ for every $B'\in \mathfrak{F}\setminus\mathfrak{F}_s$, otherwise we set $\mathfrak{F}_{s+1}:=\mathfrak{F}_s\bigcup \{B_{i_{s+1}}\}$, where $B_{i_{s+1}}$ is the ball having the largest radius among the balls non-intersecting $ \bigcup\limits_{B\in\mathfrak{F}_s}B$. Since $\mathfrak{F}$ is finite, the constructed chain stabilizes at some $s'\in\N$. We claim that $\mathfrak{F}':=\mathfrak{F}_{s'}$ is the desired subcollection. Indeed, if $\mathfrak{F}\neq\mathfrak{F}_{s'}$, then for any $B'\in \mathfrak{F}\setminus\mathfrak{F}_{s'}$, by the construction of $\mathfrak{F}_{s'}$, we have $B'\bigcap \bigcup\limits_{B\in\mathfrak{F}_{s'}}B\neq \emptyset$. Thus there exists $B\in\mathfrak{F}_{s'}$ such that $B'\bigcap B\neq \emptyset$, and hence $B'\subset 3B$.
\end{proof}

\subsection{Proofs of Theorem A and Theorem B}
\begin{proof}[Proof of Theorem A]
First, we will show for any $\alpha\in\La_X:=\{\alpha\in\R:K_{\alpha}\neq \emptyset\}$ that 
\begin{equation}\label{UB for rate function}
    I_X(\alpha)\leq h,
\end{equation}
where $h=h_{\top}(T,\Omega)$ (c.f. Lemma \ref{A-B_constant_and_top_entropy}).
This, in particular, implies that $\La_X\subset\dom(I_X)$. 

Fix any decreasing sequence $\{\delta_m\}$ of positive numbers such that $\delta_m\to 0$ as $m\to\infty$. For $\alpha\in\La_X$ set 
\begin{equation}\label{important sets}
    K_{\alpha,m}^{(n)}:=\{w\in\Omega: -\delta_m<\frac{1}{n}X_n(w)-\alpha<\delta_m\}, \text{ and } Z_{m,j}:=\bigcap_{n\geq j}K_{\alpha,m}^{(n)}\bigcap K_{\alpha}.
\end{equation}
Then for every $m\in\N$, one has $K_{\alpha}=\bigcup_j Z_{m,j}$ . 

We prove (\ref{UB for rate function}) by assuming the contradiction. 
Suppose that $I_X(\alpha)>h$ for some $\alpha\in\La_X$.
Then there is $\Delta>0$ such that $h<\Delta<I_X(\alpha)$. 
{
Take $\delta>0$, and consider the closed interval $[\alpha-\delta,\alpha+\delta]$.}
Since $I_X$ is a lower semi-continuous function, {
$I_X$ attains its minimum on $[\alpha-\delta,\alpha+\delta]$ at some point $\alpha_\delta\in [\alpha-\delta,\alpha+\delta]$, i.e., $I_X(\alpha_\delta)=I_X([\alpha-\delta,\alpha+\delta])$.
Thus for all $\delta>0$, $I_X(\alpha_\delta)\leq I_X(\alpha)$, therefore, $\limsup\limits_{\delta\to 0}I_X(\alpha_{\delta})\leq I_X(\alpha)$.
But on the other hand,  $\liminf\limits_{\delta\to 0} I_X(\alpha_\delta)\geq I_X(\alpha)$ since $I_X$ is lower semi-continuous and $\alpha_\delta\to\alpha$ as $\delta\to 0.$ 
Thus one can conclude that}
$\lim\limits_{\delta\to 0}I_X([\alpha-\delta,\alpha+\delta])=I_X(\alpha)$, therefore, there exists $\delta>0$ such that $I_X([\alpha-\delta,\alpha+\delta])> \Delta$.
Take any $\epsilon>0$, and fix it. 
Then since $\{X_n\}$ satisfies the weak Bowen condition, there exists $N_0\in\N$ such that $\frac{v_{n,\epsilon}(X_n)}{n}<\frac{\delta}{2}$ for all $n\geq N_0$. 
Furthermore, since $\delta_m\to 0$, there is $m^*\in \N$ with $\delta_{m^*}<\delta/2$. 
Finally, since $\bigcup_j Z_{j,m^*}=K_{\alpha}\neq\emptyset$ and $\{Z_{j,m^*}\}_j$ is an ascending chain of the sets, there exists $z\in K_{\alpha}$ and $j^*\in\N$ such that for all $j\geq j^*$,
$z\in Z_{j,m^*}$. 
Thus, for all $N\geq j^*$, one can readily check that
\begin{equation}\label{Bowen ball and level set}
    B_N(z,\epsilon)\subset \Big\{w\in\Omega:\frac{1}{N}X_N(w)\in\Big(\alpha-\delta_{m^*}-\frac{v_{N,\epsilon}(X_N)}{N},\alpha+\delta_{m^*}+\frac{v_{N,\epsilon}(X_N)}{N} \Big) \Big\}.
\end{equation}
Therefore, for every $N\geq \max\{j^*,N_0\}$ one has
\begin{equation}
    B_N(z,\epsilon)\subset \Big\{\frac{1}{N}X_N\in\Big[\alpha-\delta,\alpha+\delta \Big] \Big\}.
\end{equation}
Since $\nu$ is an Ahlfors-Bowen measure, for every $N\geq \max\{j^*,N_0\}$ one obtains
\begin{equation}
    \frac{1}{N}\log\nu\Big\{\frac{1}{N}X_N\in[\alpha-\delta,\alpha+\delta ] \Big\}\geq \frac{1}{N}\nu(B_N(z,\epsilon))\geq-h-\frac{\log D(\epsilon)}{N}.
\end{equation}
Then (\ref{LDUB}) leads 
\[
-h\leq \limsup_{N\to\infty}\frac{1}{N}\log\nu\Big\{\frac{1}{N}X_N\in[\alpha-\delta,\alpha+\delta ] \Big\}\leq -I_X([\alpha-\delta,\alpha+\delta ])<-\Delta.
\]
But this is a contradiction to $h<\Delta$, therefore, $I_X(\alpha)\leq h$ for every $\alpha\in\La_X$.

Now we prove the inequality (\ref{MFUB_from_LDUB}). Let $\{\delta_m\}$ and $Z_{m,j}$ be as above. Take any $\alpha\in\La_X$ and $\epsilon>0$. Since $I_X$ is a lower semi-continuous function, $\lim\limits_{\delta\to 0}I_X([a-\delta,a+\delta])=I_X(a)$ for all $a\in\R$.  Since $I_X(\alpha)<+\infty$ for $\alpha\in\La_X$, for arbitrary $\sigma>0$, there exists $m'\in\N$, and $\delta'>0$ such that, for all $m\geq m'$ and $\delta\leq\delta'$
\begin{equation}
    \Big|I_X([\alpha-\delta_m-\delta;\alpha+\delta_m+\delta])-I_X(\alpha)\Big|<\sigma.
\end{equation}
Again $\lim_{n\to\infty}\frac{v_{n,\epsilon}(X_n)}{n}=0$, choose sufficiently large $N_1\in\N$ such that $\frac{v_{n,\epsilon}(X_n)}{n}<\delta'$ for all $n\geq N_1$. From (\ref{LDUB}) one can choose $N_2\in\N$ such that for all $n\geq N_2$
\begin{equation}\label{LDUB_est}
    \nu\Big\{\frac{1}{n}X_n\in[\alpha-\delta_{m'}-\delta',\alpha+\delta_{m'}+\delta'] \Big\}\leq e^{n(-I_X([\alpha-\delta_{m'}-\delta',\alpha+\delta_{m'}+\delta'])+\sigma)}\leq e^{n(-I_X(\alpha)+2\sigma)}.
\end{equation}
Now take any $ j\in\N$ with $Z_{j,m'}\neq\emptyset$, and for all $N\geq \max\{j, N_1, N_2\}$ consider the cover $\{B_N(z,\epsilon/3)\}_{z\in Z_{m',j}}$ of $Z_{m',j}$. Clearly, it also covers $\overline{Z}_{m',j}$, therefore, by compactness, there is a finite subcover $\{B_N(z_i,\epsilon/3)\}_{i}$ of $Z_{m',j}$. Then by Vitali covering lemma (Lemma \ref{Vitali covering lemma}) there exists a disjoint subcollection $\{B_N(z_{i_k},\epsilon/3)\}_{k}$ such that 
$
Z_{m',j}\subset\cup_i B_N(z_i,\epsilon/3)\subset\cup_k B_N(z_{i_k},\epsilon).
$
Therefore, without loss of generality, one can assume that there is a finite disjoint collection $\B_N:=\{B_N(z_i,\epsilon/3)\}_{z_i\in Z_{m',j}}$ of the balls such that $Z_{m',j}\subset \bigcup_i B_N(z_i,\epsilon)$.
Since $z_i\in Z_{m',j}$, one can easily check for every $i$ that
\begin{equation}
    B_N(z_i,\epsilon)\subset \Big\{w\in\Omega:\frac{1}{N}X_N(w)\in\Big(\alpha-\delta_{m'}-\frac{v_{N,\epsilon}(X_N)}{N},\alpha+\delta_{m'}+\frac{v_{N,\epsilon}(X_N)}{N} \Big) \Big\}.
\end{equation}
Thus by the choice of $N$, for all $i$,
\begin{equation}
    B_N(z_i,\epsilon)\subset \Big\{\frac{1}{N}X_N\in(\alpha-\delta_{m'}-\delta',\alpha+\delta_{m'}+\delta' ) \Big\},
\end{equation}
and hence
\begin{equation}
    \bigcup_i B_N(z_i,\epsilon)\subset \Big\{\frac{1}{N}X_N\in(\alpha-\delta_{m'}-\delta',\alpha+\delta_{m'}+\delta' ) \Big\}.
\end{equation}
Since $\B_N:=\{B_N(z_i,\epsilon/3)\}$ is a disjoint collection of the balls, one can easily get the following from the above
\begin{equation}
    \sum_i \nu(B_N(z_i,\epsilon/3))\leq \nu\Big (\Big\{\frac{1}{N}X_N\in(\alpha-\delta_{m'}-\delta',\alpha+\delta_{m'}+\delta' ) \Big\}\Big ).
\end{equation}
Thus it follows from (\ref{LDUB_est}) that
\begin{equation}
    \sum_i \nu(B_N(z_i,\epsilon/3))\leq e^{N(-I_X(\alpha)+2\sigma)}.
\end{equation}
Since $\nu$ is an Ahlfors-Bowen measure, we can thus conclude that
\begin{equation}
    \#\B_N\leq D(\epsilon/3) e^{N(h-I_X(\alpha)+2\sigma)}.
\end{equation}
Since $\{B_N(z_i,\epsilon)\}_i$ is a cover of $Z_{m',j}$, one can immediately obtain from (\ref{non-centered def of top entropy}) that
\begin{equation}
    m(Z_{m',j},t,\epsilon,N)\leq \#\B_N\cdot e^{-Nt}.
\end{equation}
Then by combining the last two inequalities, one gets
\begin{equation}\label{important}
    m(Z_{m',j},t,\epsilon,N)\leq D(\epsilon/3) e^{-N(t-h+I_X(\alpha)-2\sigma)}
\end{equation}
Notice that $\sigma$ is not chosen yet, therefore, if $t>h-I_X(\alpha)$ then one can choose $\sigma$ in the interval $\Big(0,\frac{t-h+I_X(\alpha)}{2}\Big)$. Let $t>h-I_X(\alpha)$ be any number. Then one can obtain from (\ref{important}) that $h_{\top}(T,Z_{m',j},\epsilon)\leq t$, and notice that this inequality holds for any $\epsilon>0$ and for all large $j\in\N$. Therefore, $h_{\top}(T,Z_{m',j})\leq t$ for all sufficiently large $j$, thus by the countable stability of the topological entropy, we can get $h_{\top}(T,K_{\alpha})\leq t$ since $\{Z_{m',j}\}_j\uparrow K_{\alpha}$. Thus we can obtain the desired result since $t$ is chosen arbitrarily in $(h-I_X(\alpha),\infty)$. 
\end{proof}
\begin{proof}[Proof of Theorem B]
One can readily check that $\phi_X^*$ is an extended convex lower-semi-continuous function on $\R$. Indeed, the lower semi-continuity of $\phi_X^*$ follows from the fact that it is the \textit{superior envelope} of the family $\{u_t:t\in\R\}$ of extended continuous functions $\alpha\in\R\xrightarrow{u_t}\alpha t-\phi(t)\in\overline{\R}$. Furthermore, since $\phi_X(0)=0$, one has $\phi_X^*\geq 0$. Thus $\phi_X^*$ is a rate function in the sense of Definition \ref{ldp rate function definition}.  By Theorem \ref{Gartner-Ellis UB weak LDP}, one has the weak large deviations upper bound for the sequence $\Big\{\frac{1}{n}X_n \Big\}$ with the rate function $\phi_X^*$, i.e.,
for all compact $F\Subset\R$, one has the following
\begin{equation}
\limsup_{n\to\infty}\frac{1}{n}\log\nu\Big\{\frac{1}{n}X_n\in F \Big\}\leq -\inf_{F}(\phi_X^*)=:-\phi_X^*(F).  
    \end{equation}
    Thus one can immediately obtain the statement of the theorem from Theorem A.


\end{proof}

\section{Proof of Theorem C: Lower bound}\label{proof of lower bound-section}
The goal of this section is to prove Theorem C.
For any real number $q$ and a natural number $n$ denote
$$
\Zb_{n}(q):=\int_{\Omega}e^{qX_n}d\nu.
$$
We will employ the Large Deviations technique. Namely, for every $q$, by using Cramer's  exponential tilting method, we will construct a measure $\mu_q$ (Lemma \ref{relationship between mu_q and nu}) such that for every $\omega\in\Omega$, $n\in\N$
\begin{equation}\label{eq:upboundmeasure}
\mu_q( B_n(\omega,\epsilon)) \le \frac {e^{q X_n(\omega)}}{\Zb_n(q)} \nu( B_n(\omega,\epsilon))\times K(n,q,\epsilon),
\end{equation}
where the uniform `correction' factor $K(n,q,\epsilon)>0$ is subexponential in $n$.

Once (\ref{eq:upboundmeasure}) is established, we immediately obtain that
for all $\omega\in\Omega$, one has
$$
-\frac 1n\log \mu_q(B_n(\omega,\epsilon))
\ge -\frac {q X_n(\omega)}{n}+\frac {1}{n} \log \Zb_n(q)-\frac {1}{n} \log\nu(B_n(\omega,\epsilon))+o(1).
$$
Therefore, for $\omega\in K_\alpha$ 
we then conclude 
\begin{equation}\label{estimate lower local entropy below}
    {\underline h}_{\mu_q}(T,\omega)=\lim_{\epsilon\to 0}\liminf_{n\to\infty}
-\frac 1n\log \mu_q(B_n(\omega,\epsilon))\ge h_{\top}(T,\Omega)-q\alpha +\phi_X(q).
\end{equation}

Furthermore, for a given $\alpha$, we will identify $q$ such that $\mu_q(K_\alpha)=1$ (Lemma \ref{supported_measure_on_level_set}). As usual, $q=q_\alpha$ is such that $\alpha=\phi'(q_\alpha)$.
Therefore,
by the entropy distribution principle (Theorem \ref{entropy distribution principle}), we can conclude that
\begin{equation}\label{proof lb combining}
h_{\top}(T,K_\alpha)\ge -q_{\alpha}\alpha +\phi_X(q_\alpha)+h_{\top}(T,\Omega)= h_{\top}(T,\Omega)- \phi_X^*(\alpha) =h_{\top}(T,\Omega)-I_X(\alpha).    
\end{equation}
After combining this with the upper bound, we obtain the claim of Theorem C.
Now, let us turn to the rigorous proof.

\subsection{A fundamental inequality}
In this subsection, we shall establish a compound inequality for the dynamical systems which is introduced in Subsection \ref{expansiveness and specification}.  Recall that for any $\epsilon>0$ there exists $M_1(\epsilon)\in\N$ such that for all $ x\in\Omega$ and every $n\in\N$, $T^{n+M_1}(B_n(x,\epsilon))=\Omega$ since $T$ is strongly topologically exact transformation, and there exists $M_2(\epsilon)\in\N$ such that for all $x,x'\in\Omega$ and every $ n_1,n_2\in\N$, $B_{n_1}(x,\epsilon)\cap T^{-(n_1+M_2)}\big(B_{n_2}(x',\epsilon) \big)\neq \emptyset$ since $T$ has the weak specification property (c.f. equations (\ref{strong topological exactness}) and (\ref{weak specification condition})). Let $M_1(\epsilon)$ and $M_2(\epsilon)$ be the smallest integers satisfying these properties, respectively. For $\epsilon>0$, let 
\begin{equation}\label{the constant M}
M(\epsilon):=\max\{M_1(\epsilon),M_2(\epsilon/2)\}.
\end{equation}
One has the following lemma for the expansive dynamical systems with certain mixing properties.
\begin{lemma}\label{fundamental inequality}
For any sufficiently small $\epsilon$, there exists a positive constant $C_1(\epsilon)$ such that for any $n\in\N$ all $x\in\Omega$ and for any measurable function $Y:\Omega\to[0,+\infty)$ one has
\[
 \frac{C_1(\epsilon)^{-1}}{\nu(B_n(x, \varepsilon))}   \int_{B_n (x, \varepsilon)}    Y\circ T^{n+M} d \nu(y)\leq  \int _{\Omega}  Y d \nu \leq \frac{C_1(\epsilon)}{\nu(B_n(x, \varepsilon))}    \int_{B_n (x, \varepsilon)}    Y\circ T^{n+M}  d \nu.
\]
\end{lemma}

\begin{remark}
In the above lemma, $\epsilon$ can be chosen in the interval $(0,\rho/2)$ where $\rho$ is the expansivity constant of the transformation $T$.
\end{remark}
\begin{proof}
We assume that $T$ is expansive. Let $\rho$ denote the expansivity constant of $T$.
Consider $\epsilon\in (0,\rho/2)$, $n \in \mathbb{N} $, $ x \in \Omega$ and a measurable function $Y:\Omega\to[0,+\infty)$. 
Let $\nu_{x,n}$ be the conditional probability measure on $B_n(x,\epsilon)$, i.e., $\nu_{x,n}(A):=\nu(A \vert B_n(x,\epsilon))=\frac{\nu(A)}{\nu(B_n(x, \varepsilon))} $ for a Borel set $A\subset B_n(x,\epsilon)$. 
Set $\tau_{x,n}:=\nu_{x,n} \circ T^{-(n+M)} $, then $\tau_{x,n}$ is always a probability measure on $\Omega$ since $T^{n+M}(B_n(x,\epsilon))=\Omega$. Note that $\nu_{x,n}$ and $\tau_{x,n}$ also depend on $\epsilon$, but we omit this to avoid cumbersome notations in our subsequent expressions. It is easy to check that 
\begin{align}\label{obvious relation mu and nu}
 \int _{\Omega}    Y(z) d \tau_{x,n}(z)=&\frac{1}{\nu(B_n(x, \varepsilon))}   \int _{B_n (x, \varepsilon)}    Y (T^{n+M} z) d \nu(z).
\end{align}
Therefore, the claim of the lemma is equivalent to showing the existence of a constant $C_1(\epsilon)>0$ independent of $n$ and $x$ such that 
\begin{equation}\label{restatement of fundamental inequality}
C_1(\epsilon)^{-1}  \int_{\Omega}    Y(y) d \tau_{x,n}(y)\leq  \int_{\Omega}    Y(y) d\nu(y) \leq C_1(\epsilon)  \int_{\Omega}    Y(y) d \tau_{x,n}(y).
\end{equation}
We now present the proof of this statement in three steps. 

\textit{\underline{Step 1:}}
We start by showing that there exists a constant $C_2(\epsilon)>0$, dependent only on $\epsilon$, such that for every $k\in\N$ and $y\in\Omega$,
\begin{equation}\label{important relation mu and nu}
C_2(\epsilon)^{-1}\leq \frac{\tau_{x,n}(B_k(y,\epsilon))}{\nu(B_k(y,\epsilon))}\leq C_2(\epsilon).
\end{equation}
Consider arbitrary $k\in\N$, $y\in\Omega$. By the choice of $M:=M(\epsilon)$, we have 
\[
T^{-(n+M)} B_k(y, \varepsilon/2) \cap B_n(x, \varepsilon/2) \neq \emptyset.
\]
Then for any $z\in T^{-(n+M)} B_k(y, \varepsilon/2) \cap B_n(x, \varepsilon/2)$, one has 
$$
B_{n+M+k}(z, \varepsilon/2) \subset T^{-(n+M)} B_k(y, \varepsilon) \cap B_n(x, \varepsilon).
$$ 
Since $\nu$ is an Ahlfors-Bowen measure (c.f. (\ref{Ahlfors-Bowen-measure})), we can get
\begin{equation}\label{}
D(\epsilon/2)^{-1}e^{-(n+M+k)h}\leq\nu( B_{n+M+k}(z, \varepsilon/2) )\leq\nu( T^{-(n+M)} B_k(y, \varepsilon) \cap B_n(x, \varepsilon) ).
\end{equation}
Hence
\begin{align}\label{mu nu Bowen balls}
\begin{split}
\tau_{x,n}(B_k(y,\epsilon))=\frac{ \nu( T^{-(n+M)} B_k(y, \varepsilon) \cap B_n(x, \varepsilon) )}{ \nu (B_n(x, \varepsilon)) }&\geq \frac{D(\epsilon/2)^{-1}e^{-(n+M+k)h}}{D(\epsilon) e^{-nh}} 
\\&\geq D(\varepsilon/2)^{-1}D(\varepsilon)^{-2} e^{-M h} \nu(B_k (y, \varepsilon)).
\end{split}
\end{align}
In order to get an upper bound in (\ref{important relation mu and nu}), set $B':=T^{-(n+M)} B_k(y, \varepsilon) \cap B_n(x, \varepsilon) $. Then since $\{B_{n+M+k}(z,2\epsilon):z\in B'\}$ is a cover of $\overline{B'}$, by the the compactness of $\Omega$ and the Vitali covering lemma (Lemma \ref{Vitali covering lemma}), there exists a finite disjoint subcollection $\mathcal{G}:=\{B_{n+M+k}(z_i,2\epsilon)\}_{z_i\in B'}$ such that $\{B_{n+M+k}(z_i,6\epsilon)\}_{z_i\in B'}$ covers $B'$. Denote $\mathscr{S}:=\{T^nz_i:B_{n+M+k}(z_i,2\epsilon)\in \mathcal{G}\}$. 
Because of the disjointness of $\mathcal{G}$ and by the definition of $B'$, it is easy to check that $\sharp\mathscr{S}=\sharp\mathcal{G}$ and $\mathscr{S}$ is a $(M,2\epsilon)$-separated set. Therefore, if we denote the cardinality of a maximal $(M,2\epsilon)$-separated set in $\Omega$ by $\mathfrak{s}(M,2\epsilon)$, then it is clear that $\sharp \mathcal{G }=\sharp\mathscr{S}\leq \mathfrak{s}(M,2\epsilon)$. From these and (\ref{Ahlfors-Bowen-measure}), it follows that 
\begin{equation}
\nu(B')\leq \sum\limits_i\nu(B_{n+M+k}(z_i,6\epsilon))\leq D(6\epsilon) \mathfrak{s}(M,2\epsilon) e^{-(n+M+k)h}
\end{equation}
This bound then implies 
\begin{equation}\label{mu nu Bowen balls ub}
\tau_{x,n}(B_k(y,\epsilon))=\frac{\nu(B')}{\nu(B_n(x,\epsilon))}\leq \frac{D(6\epsilon) \mathfrak{s}(M,2\epsilon) e^{-(n+M+k)h}}{D(\epsilon)^{-1} e^{-nh}}\leq D(6\epsilon)D(\epsilon)^2\mathfrak{s}(M,2\epsilon) \nu(B_k(y,\epsilon)).
\end{equation}
Combining (\ref{mu nu Bowen balls}), (\ref{mu nu Bowen balls ub}), the constant $C_2(\epsilon):=\max\{D(\varepsilon/2)D(\varepsilon)^{2} e^{M h} ;D(6\epsilon)D(\epsilon)^2\mathfrak{s}(M,2\epsilon)\}$ satisfies the inequality (\ref{important relation mu and nu}).

\textit{\underline{Step 2:}}
In this step, we shall prove (\ref{restatement of fundamental inequality}) for continuous non-negative functions: namely, if the constant $C_1(\epsilon):=C_2(\epsilon)D(2\epsilon)D(\epsilon)$ and $\psi\in C(\Omega, \R_+)$, then 
\begin{equation}\label{fundamental ineq for continuous function}
C_1(\epsilon)^{-1}\int_{\Omega}\psi d\tau_{x,n}\leq \int_{\Omega}\psi d\nu\leq C_1(\epsilon) \int_{\Omega}\psi d\tau_{x,n}.
\end{equation}
Consider an arbitrary $N\in\N$ and let $E_N$ be a maximal $(N,2\epsilon)$-separated set in $\Omega$. Thus $\{B_N(x,2\epsilon):x\in E_N\}$ is a cover of $\Omega$ and 
$\{B_N(x,\epsilon):x\in E_N\}$ are pairwise disjoint. By successively applying (\ref{Ahlfors-Bowen-measure}) and (\ref{important relation mu and nu}) to these collections of the balls, one can obtain 
\begin{equation}
\int_{\Omega}\psi d\tau_{x,n}\geq C_2(\epsilon)^{-1} D(2\epsilon)^{-1}D(\epsilon)^{-1}\int_{\Omega}\psi d\nu-C_2(\epsilon)^{-1}(v_{N,2\epsilon}(\psi)+D(2\epsilon)D(\epsilon)v_{N,\epsilon} (\psi)).
\end{equation} 
Since $v_{N,2\epsilon}(\psi)\to 0, v_{N,\epsilon}(\psi)\to 0$ as $N\to\infty$ (c.f. Remark \ref{expansivity important property}), one obtains $\int_{\Omega}\psi d\tau_{x,n}\geq C_2(\epsilon)^{-1} D(2\epsilon)^{-1}D(\epsilon)^{-1}\int_{\Omega}\psi d\nu$. The opposite inequality is proved in exactly the same way.

\textit{\underline{Step 3:}}  
Now, let $\tilde{Y}$ be any bounded non-negative measurable function on $\Omega$. By Lusin's theorem and Tietz's extension theorem (\cite{Ercan1997}), one can approximate $\tilde{Y}$ by non-negative continuous functions $\psi_k:\Omega\to\R,\;\; k\in\N$ such that 
$\int_{\Omega}\psi_kd\nu\to \int_{\Omega}\tilde{Y}d\nu$ and $\int_{\Omega}\psi_kd\tau_{x,n}\to \int_{\Omega}\tilde{Y}d\tau_{x,n}$ as $k\to\infty$. Hence it can be  obtained from (\ref{fundamental ineq for continuous function})  that 
\begin{equation}\label{fundamental ineq for bounded function}
C_1(\epsilon)^{-1}\int_{\Omega}\tilde{Y} d\tau_{x,n}\leq \int_{\Omega}\tilde{Y} d\nu\leq C_1(\epsilon) \int_{\Omega}\tilde{Y} d\tau_{x,n}.
\end{equation}

Now we conclude (\ref{restatement of fundamental inequality}) from (\ref{fundamental ineq for bounded function}).  Consider a non-negative measurable function $Y$ on $\Omega$, and let $Y_N:=Y\cdot \mathds{1}_{Y\leq N}$ for $N\in\N$. Then $Y_N$ is a bounded non-negative function, and $\{Y_N\}$ converges monotonically to $Y$ as $N\to\infty$. Therefore, by monotone convergence theorem and from (\ref{fundamental ineq for bounded function}), one gets (\ref{restatement of fundamental inequality}).
\end{proof}

{
\subsection{A "tilted" measure and its relationship with the reference measure}
For any real number $q \in \mathbb{R}$ and a natural number $n\in\N$, define a probability measure $\mu_{n,q}$ by Cramer's tilting method as
 \[
 d \mu_{n,q} : = \frac{e^{q X_n}}{\mathbb{E}_{\nu}e^{qX_n}} d \nu=\frac{e^{q X_n}}{\Zb_n(q)} d \nu.
 \]
For any $q\in\R$, consider the sequence $\{\mu_{n,q}\}_{n\geq 1}$ of probability measures. Since $\mathcal{M}(\Omega)$ is a compact set in weak$^*$ topology, there is a convergent subsequence $\{\mu_{N_j, q}\}_j$, and we denote the weak$^{*}-$limit of this subsequence by $\mu_q$.}
In this subsection, we shall establish an important relation (c.f. Lemma \ref{relationship between mu_q and nu}) between
$\mu_q$ and $\nu$ under the conditions of Theorem C. 
First, we prove an important lemma.
\begin{lemma}\label{submultiplicativity of Z_n}
For all $n  \in \mathbb{N} $, $q \in \R$, there exists a positive constant $C_4(n,q)$ with $\lim\limits_{n \to \infty} \frac{\log C_4(n,q)} {n} =0$ such that for any $m  \in \mathbb{N} $,
\begin{align}
 C_4(n,q)^{-1}
 \leq \frac{\Zb_{n+M+m}(q)} {\Zb_n (q) \Zb_m (q)}
 \leq C_4(n,q),
\end{align}
where $M\in\N$ is as in (\ref{the constant M}).
\end{lemma}\begin{proof}
Fix a sufficiently small $\epsilon$. Consider $n,m\in\N$, $q\in\R$. 
Note that by the weak Bowen condition, each $X_k$ is a bounded function. 
In what follows, $||X_k||$ denotes the supremum norm of a function $X_k$. 
Since the sequence $X$ is almost additive, one has 
\begin{align}\label{Z submultiplicativity ub start}
  \Zb_{n+M+m} (q) & = \int_{\Omega} e^{q X_{n+M+m} (y)} d \nu(y) 
  \leq e^{|q| (A_n+A_M+||X_M||)}\int_{\Omega} e^{q X_{n} (y)}   e^{q X_{m} (T^{n+M}y)} d \nu(y).
\end{align}
Let $E$ be a maximal $(n, \varepsilon)$-separated set in $\Omega$.
Then $\bigcup\limits_{x \in E} B_{n}(x, \varepsilon) = \Omega$ and $\{B_n(x,\epsilon/2)\}_{x\in E}$ is a disjoint collection of Bowen balls.
Therefore,
\begin{align}
\begin{split}
\int_\Omega e^{q X_n(y)} e^{q X_m (T^{n+M}y)} d\nu(y)
& \leq \sum_{x \in E} \int_{B_n(x, \varepsilon)}  e^{q X_n(y)} e^{q X_m (T^{n+M}y)} d\nu(y)\\
  & \leq  e^{|q| v_{n,\epsilon}(X_n) } \sum_{x\in E}  e^{q X_n(x)} \int_{B_{n}(x, \varepsilon) }  e^{q X_{m} (T^{n+M}y)} d \nu(y) \\
   &\label{Z submultiplicativity ub integral} \leq e^{|q| v_{n,\epsilon}(X_n)}C_1( \varepsilon) \Zb_m(q) \sum_{x\in E}  e^{q X_n(x)}  \nu(B_{n}(x, \varepsilon)),
\end{split}
\end{align}
  where the second inequality holds since $\{X_n\}$ satisfies weak Bowen's condition and the last inequality holds by Lemma \ref{fundamental inequality} applied to $Y=X_m$.
  Now we estimate the sum $\sum\limits_{x\in E}  e^{q X_n(x)}  \nu(B_{n}(x, \varepsilon))$ from above.
 One has
 \begin{align}
   \begin{split}
   \sum_{x\in E} e^{q X_n(x)} \nu(B_{n}(x, \varepsilon))& \leq \sum_{x\in E}  e^{q X_n(x)} D(\varepsilon/2) D(\varepsilon) \nu(B_n(x, \varepsilon/2))\\
   & \leq D(\varepsilon/2) D(\varepsilon) e^{|q| v_{n,\epsilon/2}(X_n)} \sum_{x\in E} \int_{B_n(x, \varepsilon/2) } e^{q X_{n} (y)}  d \nu \\
   &\label{Z submultiplicativity ub sum} \leq   D(\varepsilon/2) D(\varepsilon) e^{ |q| v_{n,\epsilon/2}(X_n) }  \Zb_n(q).
   \end{split}
\end{align}
After combining (\ref{Z submultiplicativity ub start}), (\ref{Z submultiplicativity ub integral}) and (\ref{Z submultiplicativity ub sum})
\begin{equation}\label{Z submultiplicativity ub}
\frac{\Zb_{n+M+m} (q) }{\Zb_n(q) \Zb_m(q)} \leq  C_5(n,q),
\end{equation}
where $C_5(n,q):=e^{|q| (A_n+A_M+||X_M||)} e^{ |q| (v_{n,\epsilon/2}(X_n)+v_{n,\epsilon}(X_n))} C_1(\varepsilon)   D(\varepsilon/2) D(\varepsilon)$.

With similar reasoning, one gets the following inequalities
\begin{align}
  \Zb_{n+M+m} (q) 
   & \geq e^{-|q| (A_n+A_M+||X_M||)} \int_{\Omega} e^{q X_{n} (y)}   e^{q X_{m} (T^{n+M}y)} d \nu(y),
\end{align}
\begin{align}
   \int_{\Omega} e^{q X_{n} (y)}   e^{q X_{m} (T^{n+M}y)} d \nu(y) 
   & \geq  e^{-|q| v_{n,\epsilon/2}(X_n) } C_1(\varepsilon/2)^{-1} \Zb_m(q) \sum_{x\in E}  e^{q X_n(x)}  \nu(B_{n}(x, \varepsilon/2)),
   \end{align}
 \begin{align}
   \sum_{x\in E}  e^{q X_n(x)} \nu(B_{n}(x, \varepsilon/2))
   & \geq e^{ - |q| v_{n,\epsilon}(X_n) }   D(\varepsilon/2)^{-1} D(\varepsilon)^{-1} \Zb_n(q),
\end{align}
and after combining these inequalities 
\begin{equation}\label{Z submultiplicativity lb}
\frac{\Zb_{n+M+m} (q) }{\Zb_n(q) \Zb_m(q)} \geq C_6(n,q)^{-1},
\end{equation}
where $ C_6(n,q):=e^{|q| (A_n+A_M+||X_M||)} e^{ |q|(v_{n,\epsilon/2}(X_n)+v_{n,\epsilon}(X_n))} C_1 (\varepsilon/2)  D(\varepsilon/2)  D(\varepsilon)$.

From (\ref{Z submultiplicativity ub}) and (\ref{Z submultiplicativity lb}), it is clear that $C_4(n,q):=\max\{C_5(n,q),C_6(n,q)\}$ satisfies the conditions of the lemma.
\end{proof}

\begin{lemma}\label{relationship between mu_q and nu}
For any sufficiently small $\epsilon>0$ and for any $n\in\N$, $q\in\R$ there exists a constant $C_7(n,q,\epsilon)>0$ such that $\lim\limits_{n\to\infty}\frac{\log C_7(n,q,\epsilon)}{n}=0$ and for all $x\in\Omega$,
\begin{equation}\label{the equation in lemma_relationship between mu_q and nu}
  \mu_q(B_n(x, \varepsilon)) \leq  C_7(n,q,\epsilon)  \frac{e^{qX_n(x)}}{\Zb_n(q)} \nu(B_n(x, \varepsilon)).
\end{equation}
\end{lemma}
\begin{proof}
Consider again $\epsilon\in(0,\rho/2)$, where $\rho$ is the expansivity constant of the transformation $T$, and $x\in\Omega$, $n\in\N$, $q\in\R$. 
By the Portmanteau theorem, we have
\begin{align}\label{Portmanteau Theorem}
  \mu_q(B_n (x,\varepsilon)) \leq 
  \liminf_{j \to \infty}  \int_{B_n(x, \varepsilon)} \frac{e ^{q X_{N_j}(y)}}{\Zb_{N_j}(q)} d \nu(y)
  \end{align}
 Now we estimate the numerator $\int_{B_n(x, \varepsilon)} e ^{q X_{N_j}(y)}d \nu(y)$ and denominator $\Zb_{N_j}(q)$ separately.
By Lemma \ref{submultiplicativity of Z_n}, one has 
\begin{equation}\label{applying submultiplicativity of Z_n}
\frac{1}{\Zb_{n+M+N_j-n-M}(q)} \leq C_4(n, q) \frac{1}{\Zb_{n}(q)  \Zb_{N_j-n-M}(q)}.
\end{equation}

For the integral in the numerator, the following estimations are valid
  \begin{align}\label{relationship mu_q nu integral numerator starting}
  \begin{split}
  \int_{B_n(x, \varepsilon)} e ^{q X_{n+M+N_j-n-M}(y)} d \nu(y) &\leq
  e^{|q|(A_n+A_M)} \int_{B_n(x, \varepsilon)} e ^{q X_{n}(y)}  e ^{q X_{M}(T^n y)} e ^{q X_{N_j-n-M}(T^{n+M}y)}d \nu(y)\\
  & \leq e^{ |q|(A_n+A_M+||X_M||)}  \int_{B_n(x, \varepsilon)} e ^{q X_{n}(y)} e ^{q X_{N_j-n-M}(T^{n+M}y)}d \nu(y)\\
  & \leq  C_8(n,q,\epsilon)e^{q X_n (x)}  \int_{B_n(x, \varepsilon)} e ^{q X_{N_j-n-M}(T^{n+M}y)}d \nu(y),
  \end{split}
 \end{align}
 where $C_8(n,q,\epsilon):=e^{|q|(A_n+A_M+||X_M||+ v_{n,\epsilon}(X_n)) }$.
 Note that 
 \begin{equation}
 \int_{B_n(x, \varepsilon)} e ^{q X_{N_j-n-M}(T^{n+M}y)}d \nu(y)\leq C_1(\varepsilon)   \Zb_{N_j-n-M}(q) \nu(B_n(x, \varepsilon)).
 \end{equation}
 Thus  it follows from (\ref{relationship mu_q nu integral numerator starting})  that
 \begin{equation}
 \int_{B_n(x, \varepsilon)} e ^{q X_{N_j}(y)} d \nu(y) \leq C_1(\epsilon)C_8(n,q,\epsilon)e^{q X_n (x)}\nu(B_n(x, \varepsilon))\Zb_{N_j-n-M}(q).
 \end{equation}
 Hence it easy to see from (\ref{Portmanteau Theorem})-(\ref{applying submultiplicativity of Z_n}) that the constant $C_7(n,q,\epsilon):=C_1(\epsilon)C_4(n,q)C_8(n,q,\epsilon)$ satisfies (\ref{the equation in lemma_relationship between mu_q and nu}), and by Lemma \ref{submultiplicativity of Z_n} and the weak Bowen condition, $\lim\limits_{n\to\infty}\frac{\log C_7(n,q,\epsilon)}{n}=0$.
 \end{proof}

\subsection{Differentiability of the log-moment generating function}
\begin{lemma}\label{existence of moment-generating function}
 For all $q\in\R$,  the limit $\lim\limits_{n\to\infty} \frac{1}{n}\log\Zb_n(q)$ exists and finite. Thus, 
 $\phi_X=\lim\limits_{n\to\infty} \frac{1}{n}\log\Zb_n$ is a finite convex function. Furthermore, $\phi_X=I_X^*$ and $\phi_X^*=I_X$.
 \end{lemma}

First, we state a slight generalization of the classical Fekete's lemma which is used to verify the existence and finiteness of the limit $\lim\limits_{n\to\infty} \frac{1}{n}\log\Zb_n(q)$. 

\begin{lemma}\cite{Fan1996}\label{Fekete's lemma}
    Let $\{a_n\}_{n\geq 0}$ be a sequence of real numbers. Assume that there exists a natural number $N\in\N$ and real numbers $D$, $\Delta_n\geq 0$, $n\in\N$ with $\lim_{n\to\infty} \frac{\Delta_n}{n}=0$ such that for all $m\in \Z_+$ and $n\geq N$
\begin{equation}\label{weak form of subadditivity}
    a_{n+m}\leq a_n+a_m+D+\Delta_n.
\end{equation}
Then the limit $\lim_{n\to\infty}\frac{a_n}{n}$ exists and belongs to $\R\bigcup\{-\infty\}$.
\end{lemma}


\begin{proof}[Proof of Lemma \ref{existence of moment-generating function}]
By using Lemma \ref{submultiplicativity of Z_n}, it can be shown for $q\in\R$, $n\geq M+1$ and any $m\in\Z_+$ that
\begin{equation}\label{almost additivity log Zb}
    \log\Zb_{n+m}(q)\leq \log\Zb_n(q)+\log\Zb_m(q)+2\log C_4(n-M,q).
\end{equation}
Thus by Lemma \ref{Fekete's lemma}, the limit $\lim\limits_{n\to\infty} \frac{1}{n}\log\Zb_n(q)$ exists and $\lim\limits_{n\to\infty} \frac{1}{n}\log\Zb_n(q)<+\infty$. 
By applying a similar argument to the sequence $\{-\log\Zb_n(q)\}$, one shows that $\lim\limits_{n\to\infty} \frac{1}{n}\log\Zb_n(q)>-\infty$.

Clearly, each $\log\Zb_n$ is a convex function, thus $\phi_X$ is also a convex function as the pointwise limit of convex functions.

Since for all $q\in\R$, $\phi_X(q)<+\infty$, and $\{\frac{1}{n}X_n\}_{n\in\N}$ satisfies LDP with the good rate function $I_X$, we immediately obtain from Theorem \ref{Varadhan theorem} that $I_X^*=\phi_X$. Note that $I_X=\phi_X^*$ also holds by the Fenchel duality lemma (\cite{Dembo-Zeitouni-Book}), since $I_X$ is convex and lower semi-continuous. 
 \end{proof}

Now we show that $\phi_X$ is a differentiable function. In fact, the following holds true.
\begin{lemma}\label{differentiability of moment-generating function}
The function $\phi_X$ is continuously differentiable on $\R$, i.e., $\phi_X\in C^1(\R)$.
\end{lemma}
\begin{proof}
By Lemma \ref{existence of moment-generating function}, we have that $I_X^*=\phi_X$. Thus since $I_X$ is an essentially strictly convex rate function, $\phi_X$ is an essentially smooth function by Theorem \ref{the theorem in Rockafellar book}. Hence since $\dom(\phi_X)=\R$ by Lemma \ref{existence of moment-generating function}, we have that $\phi_X$ is a differentiable function on $\R$. Note that $\phi_X$ is a convex function, therefore, $\phi_X\in C^1(\R)$ by Darboux's intermediate value theorem.  
\end{proof}
By the monotonicity of $\phi_X'$, we define $\underline{\alpha}:=\lim_{t\to -\infty}\phi_X'(t)$ and $\overline{\alpha}:=\lim_{t\to +\infty}\phi_X'(t)$. Then by convexity of $\phi_X$,  it is easy to verify the following limits
\begin{equation}\label{L'Hoptials rule at infinity}
\underline{\alpha}=\lim\limits_{q\to-\infty}\frac{\phi_X(q)}{q}, \text{ and } \overline{\alpha}=\lim\limits_{q\to+\infty}\frac{\phi_X(q)}{q}.
\end{equation}




\subsection{
{A 
measure supported on the level set and a proof of Theorem C}}

\begin{lemma}\label{supported_measure_on_level_set}
Under the conditions of Theorem C, for every $\alpha\in (\underline{\alpha}, \overline{\alpha})$, there exists $q\in\R$ such that $\mu_q(K_{\alpha})=1$. In particular, $(\underline{\alpha}, \overline{\alpha})\subset \L_X$.
\end{lemma}

The proof of the lemma relies on the following well-known result in the theory of convex functions.

\begin{prop}[\cites{Ellis-Book, Rockafellar-Book}]\label{Differentiable}
For any convex function $\Lambda$ defined on $\mathbb{R}$, we have \\
$(1)$ $\beta q \leq \Lambda^{\ast} (\beta) + \Lambda(q)  \   \ \text{for any  } \; q, \beta \in \mathbb{R}$,\\
$(2)$ if $\Lambda$ is differentiable at $q$ then $\beta q = \Lambda^{*} (\beta) + \Lambda(q)  \Longleftrightarrow   \beta = \Lambda' (q)$.
\end{prop}
\begin{proof}[ Proof of Lemma \ref{supported_measure_on_level_set}]
Note that the for every $q\in\R$, the limit $\lim\limits_{n \to \infty} \frac{1}{n} \log \mathbb E_{\nu} e^{qX_n}$ exists by Lemma \ref{existence of moment-generating function}. One (\cite{Fan2021}) can prove that the following limit
\[
C_{q} (\beta) : = \lim_{n \to \infty} \frac{1}{n} \log \mathbb E_{\mu_q} e^{\beta X_n},
\]
also exists, and $C_{q} (\beta)= \phi_X(q+\beta) -\phi_X(q), \; \forall \beta,q \in\R$.
Then we have $C_{q}'(\beta) = \phi_X'(q+\beta) $ and $C_{q} (0)=0$, and thus $C_{q} ^{*} (\gamma) \geq 0$ for any $\gamma \in \mathbb{R}$.
Moreover, $C_{q}^{\ast} (\alpha) = 0$ if and only if $\alpha= C_{q}'(0) = \phi_X' (q)$.
The following two observations show that $\mu_q$ is supported on the level set $K_{\alpha}$:
\begin{itemize}
\item For any nonempty closed set $F\subset\R$ such that $\alpha\notin F$ , we have $\eta_F:=\inf\limits_{\gamma \in F} C_{q}^{*} (\gamma) >0$.

\item Since the limit $\lim_{n \to \infty} \frac{1}{n} \log \mathbb E_{\mu_q} e^{\beta X_n}$ exists for any $\beta \in \mathbb{R}$, by the first part of Theorem \ref{Gartner-Ellis}, for any closed set $F$,  we have,
\[
\limsup_{n \to \infty} \frac{1}{n} \log \mu_q (\{ \omega \in \Omega : \frac{1}{n} X_n (\omega) \in F \} ) \leq -\eta_F<0.
\]
\end{itemize}
In fact, for $m\in\N $, consider $F_m:=\R\setminus (\alpha-\frac{1}{m},\alpha+\frac{1}{m})$. From the above observations and the Borel-Cantelli lemma, one has for every $m\in\N$ that
\begin{equation}
    \mu_q\Big(\bigcap_{j\in\N}\bigcup_{n\geq j}\Big\{\frac{1}{n}X_n\in F_m \Big\}\Big)=0.
\end{equation}
Therefore,
\begin{equation}
\mu_q(K_{\alpha})=\lim\limits_{m\to\infty}\mu_q\Big(\bigcup_{j\in\N}\bigcap_{n\geq j}\Big\{\frac{1}{n}X_n\in F_m^c \Big\}\Big)=1.
\end{equation}

\end{proof}

\begin{proof}[Proof of Theorem C] 
%
\textit{\underline{The first claim of Theorem C:}}  $(\underline{\alpha},\overline{\alpha})\subset \L_X$ is proven in Lemma \ref{supported_measure_on_level_set}. Now we show that $\L_X\subset [\underline{\alpha},\overline{\alpha}]$. Take any $\alpha\in\L_X$, then by Theorem B,  $\phi_X^{*}(\alpha)\leq h=h_{\top}(T,\Omega)$. Thus by Young's inequality (the first statement of Proposition \ref{Differentiable}), one has for all $q\in\R$,
\begin{equation}
\alpha q\leq \phi_X^*(\alpha)+\phi_X(q)\leq h+\phi_X(q).
\end{equation}
Hence for $q<0$,
\begin{equation}\label{L_X interval}
\alpha\geq \frac{h}{q}+\frac{\phi_X(q)}{q}, \text{ and for } q>0,\; \; \alpha\leq \frac{h}{q}+\frac{\phi_X(q)}{q}.
\end{equation}
Furthermore, from (\ref{L'Hoptials rule at infinity}) and (\ref{L_X interval}), for $q<0$, we conclude that $\alpha\geq \lim\limits_{q\to-\infty}\big(\frac{h}{q}+\frac{\phi_X(q)}{q}\big)=\underline{\alpha}$, and for $q>0$, we conclude that $\alpha\leq \lim\limits_{q\to+\infty}\big(\frac{h}{q}+\frac{\phi_X(q)}{q}\big)=\overline{\alpha}$.


\textit{\underline{The second claim of Theorem C:}} By combining Lemma \ref{relationship between mu_q and nu}, Lemma \ref{supported_measure_on_level_set}, Proposition \ref{Differentiable} and Theorem \ref{entropy distribution principle}, one concludes the proof of the second claim (c.f. (\ref{eq:upboundmeasure})-(\ref{proof lb combining})).



\end{proof}


\section{Examples and Discussion}\label{examples}
In this section, we shall discuss some of the known results of multifractal formalism, and show how they can be treated by the results of this paper.

\subsection{Lyapunov exponents for products of matrices.}
We shall adopt the setup in (\cite{Feng2003},\cite{Feng2009}, \cite{Feng-Lau2002}). Consider a map $M:\Omega\to M_d(\R)$, where $\Omega:=\{1,...,q\}^{\Z_+}$ is a full shift space with the shift transformation $T$, and $M_d(\R)$ is the algebra of $d\times d$ real matrices. For $n\in\N$, $x\in\Omega$, set $M_n(x):=M(T^{n-1}x)\;M(T^{n-2}x)\cdots M(x)$. 
We impose the following two conditions on the map $M$:
\begin{itemize}
    \item[(C1):] $M_n(x)\neq O$ holds for all $n\in\N$ and $x\in\Omega$, where $O$ is the zero matrix;
    \item[(C2):] $M$ depends only on finitely many coordinates, i.e., $M$ is locally constant.
\end{itemize}
Note that condition (C1) is automatically satisfied if, for example, the map $M$ takes values in the set $GL_d(\R)$ of irreducible matrices, or the set $M_d^+(\R)$ of $d\times d$ positive matrices. In \cite{Feng2009}, the authors imposed an irreducibility condition which also implies the condition (C1).
Under the condition (C1), one can define the (upper) Lyapunov exponent of the cocycle $(T, M)$ as: for $x\in\Omega$,
$$
\lambda^+_M(x):=\lim_{n\to\infty}\frac{1}{n}\log ||M_n(x)||,
$$
provided the limit exists, where $||\cdot||$ is a matrix norm.
Let $K_{\alpha}:=\{x\in\Omega: \lambda^+_M(x)=\alpha\}$, and $\La_M:=\{\alpha\in\R: K_{\alpha}\neq \emptyset\}$.
We denote $X_n(x):=\log ||M_n(x)||$ for all $n\in\N$, $x\in\Omega$. Then it is clear that $X_n$ is a locally constant function since $M$ is only dependent on finitely many coordinates, therefore, the sequence $\{X_n\}$ satisfies the Bowen condition.
In \cite{Feng-Lau2002}, the authors introduced the pressure function $P_M$ of the cocycle $(T,M)$ by
$$
P_M(t):=\lim_{n\to\infty}\frac{1}{n}\log\sum_{\omega\in\Omega_n}\sup_{x\in[\omega]}||M(T^{n-1}x)\;M(T^{n-1}x) \cdots M(x)||^t, \ \ \forall~t\in\R,
$$
where $\Omega_n:=\{1,...,q\}^n$, $n\in\N$. Note that the above limit exists and 
\begin{equation}\label{lyapunov exp pressure and mom generate funct}
    P_M(t)=\log q+\phi_X(t),\; t\in\R,
\end{equation}
where $\phi_X$ is the log moment generating function of the sequence $\{X_n\}$, i.e., $\phi_X(t)=\lim_{n\to\infty}\frac{1}{n}\log\int_{\Omega}\exp{(tX_n)}d\nu,\;\; t\in\R$, and $\nu$ is the uniform Bernoulli measure on $\Omega$. Then by applying Theorem B, one can obtain that 
\begin{equation}\label{lyapunov exp_ub}
    h_{\top}(T,K_{\alpha})\leq \log q-\phi_X^*(\alpha)=\inf\limits_{q\in\R}\{-\alpha q+P_M(q)\}
\end{equation}
for all $\alpha\in\La_M$. Note that the last inequality generalizes the upper bound of Theorem 1.1 in \cite{Feng2009}. In the same theorem, the authors also proved the lower bound under much stronger conditions than the conditions we have imposed above. In this setup, we can not treat the lower bound with Theorem C, since the sequence $\{X_n\}$ is not almost additive in general.  However, if the map $M$ takes values in the set $M_d^+(\R)$ of $d\times d$ positive matrices, then it was proven in \cite{Feng-Lau2002}, \cite{Feng2003} that the sequence $\{X_n\}$ is in fact almost additive. Furthermore, in this case, D.J. Feng has proven in \cite{Feng2003} that the pressure function $P_M$ is differentiable, thus $\phi_X$ is also differentiable by (\ref{lyapunov exp pressure and mom generate funct}). Then, in this particular setup, one can apply Theorem C and obtain the equality in (\ref{lyapunov exp_ub}).



\subsection{Ratios of the Birkhoff ergodic averages}

In large deviation theory, the contraction principle describes how the large deviation rate function evolves if the corresponding sequence of random elements (probability measures, random variables) is transformed by a continuous map. More specifically, if $\mathcal{X}$ and $\mathcal{Y}$ are Hausdorff topological spaces and a family $\{\mu_n\}$ of probability measures on $\mathcal{X}$ satisfies the LDP with the rate function $I$, $F:\mathcal{X}\to \mathcal{Y}$ is a continuous map, then $\tau_n:=\mu_n\circ F^{-1}$ satisfies LDP with the rate function $J(y):=\inf_{x:\;F(x)=y} I(x)$. For example, suppose $X_n^{(1)}=S_n\phi$, $X_n^{(2)}=S_n\psi$, where $\phi,\psi$ are continuous functions on a compact metric space $\Omega$ and $\psi$ is strictly positive. Suppose $\nu$ is a probability measure such that the (vector) sequence $X=(X_n^{(1)}, X_n^{(2)})$ satisfies the LDP with the rate function $I_X(\alpha,\beta)$. Then a new sequence of random variables 
$$
Z_n:=\frac{X_n^{(1)}}{X_n^{(2)}}=\frac{S_n\phi}{S_n\psi},\;\; n\geq 1,
$$
satisfies LDP with the rate function $J_Z(\gamma)=\inf_{\alpha/\beta=\gamma}I_X(\alpha,\beta)$. 
Thus we can easily obtain many non-trivial observables satisfying LDP. In case, $X_n^{(1)}=S_n\phi$, $X_n^{(2)}=S_n\psi$, $\psi>0$, one easily shows that $\{nZ_n\}$ satisfies {the weak Bowen condition}, therefore, Theorem A applies, and hence $h_{\top}(T,K_{\gamma}^{Z})\leq h_{\top}(T,\Omega)-\inf_{\alpha/\beta=\gamma}I_X(\alpha,\beta)=h_{\top}(T,\Omega)-J_Z(\gamma)$. 
In this case, one can also easily obtain a lower bound as well. 
Indeed, for any $\alpha,\beta$ such that 
$$
K_{\alpha,\beta}^X:=\{\omega\in\Omega: \frac{1}{n}(S_n\phi)(\omega)\to\alpha,\; \frac{1}{n}(S_n\psi)(\omega)\to\beta\},
$$ 
we have that $K^X_{\alpha,\beta}\subseteq K^Z_{\alpha/\beta}$. 
Therefore, $h_{\top}(T,K^Z_{\gamma})\geq h_{\top} (T,K^X_{\alpha,\beta})$, where $\gamma=\alpha/\beta$.
Note that Theorem C has been established for random variables, but the same argument works for random vectors as well. (There are also other results for the level sets of the form $K^X_{\alpha,\beta}$, c.f., \cite{Takens-Verbitskiy2003}.) Hence $h_{\top}(T,K_{\gamma}^Z)\geq h_{\top}(T,\Omega)-I_X(\alpha,\beta),$ and thus 
\begin{align*}
   h_{\top}(T,K_{\gamma}^Z)&\geq \sup_{\alpha/\beta =\gamma}[h_{\top}(T,\Omega)-I_X(\alpha,\beta)]\\&= h_{\top}(T,\Omega)-\inf_{\alpha/\beta=\gamma}I_X(\alpha,\beta)\\&=h_{\top}(T,\Omega)- J_Z(\gamma).
\end{align*}
Therefore, combining two inequalities conclude that $h_{\top}(T,K_{\gamma}^Z)=h_{\top}(T,\Omega)- J_Z(\gamma)$.

\subsection{Almost additive sequences over a subshift of finite type.}
In \cite{Bomfim-Varandas2015}, the authors considered almost additive sequences of continuous potentials and they obtained an expression for the multifractal spectra of the level sets of the sequence in terms of a rate function. In this subsection, we shall discuss that the setup considered in \cite{Bomfim-Varandas2015} is covered by Theorem C. 
Let us now recall the result of \cite{Bomfim-Varandas2015} using the notation of the present paper.
\begin{theorem}[\cite{Bomfim-Varandas2015}]\label{Bomfim-Varandas2015}
Let $(\Omega, T)$ be a topologically mixing subshift of finite type. Assume that $X:=\{X_n\}_{n\in\N}$ is an almost additive sequence satisfying Bowen's condition, and $X$ is not cohomologous to a constant, i.e., $\frac{X_n}{n}$ does not uniformly converge  to a constant as $n\to\infty$. Furthermore, assume that $\lim_{n\to\infty}\frac{1}{n}\int_{\Omega}X_nd\nu=\inf_{n\in\N}\frac{1}{n}\int_{\Omega}X_nd\nu=0$, where $\nu$ is the unique measure of maximal entropy for $T$. Then
$$h_{\top}(T,K_{\alpha})=h_{\top}(T,\Omega)-\phi_X^*(\alpha).$$
\end{theorem}
It is clear that the Bowen condition implies the weak Bowen condition, and a mixing subshift of finite type satisfies the condition (A3) in Theorem C. In Theorem C, it is also assumed that the sequence $\{\frac{1}{n}X_n\}_{n\in\N}$ satisfies LDP with an essentially strictly convex good rate function. Note that if the moment generating function $\phi_X$ of the sequence $X$ is finite and differentiable on $\R$, then by Theorem \ref{Gartner-Ellis} $\{\frac{1}{n}X_n\}_{n\in\N}$ satisfies LDP with a good rate function $\phi_X^*$. Furthermore, by Theorem \ref{the theorem in Rockafellar book}, $\phi_X^*$ is also an essentially strictly convex function. 
Below, we show that the conditions of Theorem \ref{Bomfim-Varandas2015} imply differentiability of the moment-generating function $\phi_X$.
In \cite{Bomfim-Varandas2015}, it is observed that 
\begin{equation}\label{relation_mom_gen_func_and_pressure_func_almost_additive_setup}
\phi_X(t)=P_{top}(T,tX)-h_{\top}(T,\Omega), \text{ for all }t\in\R,
\end{equation}
where $P_{top}(T,tX)$ is \textit{the topological pressure} (for the concept, see  \cite{Bomfim-Varandas2015}, \cite{Barreira2006}, \cite{Feng-Huang2010}, \cite{Zhao-Zhang-Cao2011} ) \textit{of $T$ of the almost additive sequence $X=\{X_n\}$}, and $\phi_X(t)=\lim\limits_{n\to\infty}\frac{1}{n}\log\int\limits_{\Omega}e^{tX_n(x)}\nu(dx)$ as usual.  In \cite{Barreira-Doutor2009}, by generalising the well-known result in \cite{Takens-Verbitskiy1999} to almost additive setting, the authors proved that $t \mapsto P_{top}(T,tX)$ is differentiable at $0$ if the almost additive sequence $X$ has unique \textit{equilibrium state} (for the concept, see \cite{Bomfim-Varandas2015}, \cite{Barreira2006}, \cite{Feng-Huang2010}, \cite{Zhao-Zhang-Cao2011}), and showed
$$
\frac{d}{dt}P_{top}(T,tX)\big\rvert_{t=t_0}=\lim_{n\to\infty}\int_{\Omega}\frac{X_n}{n}d\mu_{t_0X},
$$
where $\mu_{t_0X}$ is the unique  equilibrium state of $t_0X=\{t_0X_n\}$.
In 2006, L. Barreira proved that if an almost additive sequence $X$ satisfies Bowen's condition, then it has a unique equilibrium state \cite{Barreira2006}, therefore, from the above one can conclude that the pressure function $t\in\R \mapsto P_{top}(T,tX)$ is differentiable everywhere on $\R$. Then from (\ref{relation_mom_gen_func_and_pressure_func_almost_additive_setup}), the log-moment generating function $\phi_X$ is also differentiable on $\R$, and hence, the conditions of Theorem \ref{Bomfim-Varandas2015} imply the conditions of Theorem C.

\subsection{Multifractal Formalism in the absence of strict convexity of the rate function.} 
The Manneville--Pomeau map is one of the well-known examples of a non-uniformly expanding interval map. For a given $s\in(0,1)$, the Manneville--Pomeau map $T:[0,1]\to [0,1]$ is defined as 
\[
T(x):=x+x^{1+s} \mod 1.
\]
It should be noted that $T$ is topologically conjugate to a one-sided full shift, therefore, it has the specification property. It is also easy to check that $T$ is strongly topologically exact and expansive. The \textit{upper Lyapunov exponent} $\lambda_{T}$ is defined as 
\[
\lambda_T(x)=\lim\limits_{n\to\infty}\frac{1}{n}\log (T^n)'(x)
\]
provided that the limit exists at a point $x\in[0,1]$. Set 
\[
K_{\alpha}:=\{x\in\Omega: \lambda_T(x)=\alpha\},\;\;\;\; \La_T:=\{\alpha\in\R: K_{\alpha}\neq \emptyset\}.
\]
Define, $\psi(x):=\log T'(x)=\log(1+(1+s)x^s), \; x\in[0,1]$. Then it is easy to notice that the Lyapunov exponent $\lambda_T$ exists at a point $x$ if and only if the limit $\lim\limits_{n\to\infty}\frac{1}{n}S_n\psi(x)$ converges, and $\lim\limits_{n\to\infty}\frac{1}{n}S_n\psi(x)=\lambda_T(x)$. It was shown in \cite{Prellberg-Slawny1992}, \cite{Urbanski1996} that the pressure function $P_{\psi}(q):=P(T,q \psi), \; q\in\R$ exhibits a \textit{first order phase transition} at the point $q=-1$, in a sense that $P_{\psi}$ is not differentiable at $q=-1$, furthermore, it is strictly convex positive differentiable function in $(-1,\infty)$, and  $P_{\psi}(q)=0$ for all $q<-1$. Since $T$ is expansive and satisfies the specification, there is a unique measure of maximal entropy $\nu\in\M([0,1])$ for $T$, and it also becomes the Ahlfors-Bowen measure (\cite{Walters2000}). 
Define $X_n:=S_n\psi$ for all $n\in\N$, and set $\phi_X(q):=\limsup\limits_{n\to\infty}\frac{1}{n}\log \int_{0}^{1} e^{qX_n}d\nu,\; q\in\R$. Then one can easily check that $\phi_X(q)=h_{\top}(T,[0,1])+P_{\psi}(q)$ for all $q\in \R$.  Thus since
 $\{X_n\}$ satisfies the weak Bowen condition one obtains from Theorem B that for all $\alpha\in\La_{T}$,
\begin{equation}\label{manneville-pomeau-ub}
    h_{\top}(T, K_{\alpha})\leq h_{\top}(T,[0,1])-\phi_X^*(\alpha)=-P_{\psi}^*(\alpha)=\inf\limits_{q\in\R}\{P(T,q\psi)-\alpha q\}.
\end{equation}
 Now we discuss the lower bound. As have we already mentioned, $\phi_X$ is not differentiable everywhere; otherwise, the sequence $\big\{\frac{1}{n}X_n\big\}$ satisfies LDP with an essentially strictly convex good rate function. However, Lemma \ref{relationship between mu_q and nu} and Lemma \ref{supported_measure_on_level_set} are still applicable if $\alpha\in \{0\}\cup ((P_{\psi})'_{+}(-1),\;\;\overline{\alpha})$, where $(P_{\psi})'_{+}(-1)>0$ is the right derivative of $P_{\psi}$ at the point $-1$, and $\overline{\alpha}:=\lim\limits_{q\to\infty}P'_{\psi}(q)$. Thus one can obtain from these lemmas that 
 \begin{equation}\label{manneville-pomeau-lb}
    h_{\top}(T, K_{\alpha})= h_{\top}(T,[0,1])-\phi_X^*(\alpha)=-P_{\psi}^*(\alpha)=\inf\limits_{q\in\R}\{P(T,q\psi)-\alpha q\},
\end{equation}
holds for all $\alpha\in \{0\}\cup ((P_{\psi})'_{+}(-1),\;\;\overline{\alpha})$. 

We should note that in \cite{Takens-Verbitskiy2003}, the above equality was proven for all $\alpha\in (0, \overline{\alpha})$ by a different method.


\subsection{Local or pointwise entropies.}
Recall that the upper and the lower local entropies of an invariant measure $\mu$ at a point $x\in\Omega$ are defined as 
$$
\overline{h}_{\mu}(T,x):=\lim_{\epsilon\to 0+}\limsup_{n\to\infty}-\frac{1}{n}\log \mu(B_n(x,\epsilon)),\;\; \underline{h}_{\mu}(T,x):=\lim_{\epsilon\to 0+}\liminf_{n\to\infty} -\frac{1}{n}\log \mu(B_n(x,\epsilon)).
$$
If $\overline{h}_{\mu}(T,x)=\underline{h}_{\mu}(T,x)$, we say that the local entropy of the measure $\mu$ at $x$ exists and denote by $h_{\mu}(T,x)$ the common value. Finally, consider the $\alpha$-level set $K^{(\mu)}_{\alpha}$ of the local entropies 
$$
K^{(\mu)}_{\alpha}:=\{x\in\Omega:h_{\mu}(T,x)=\alpha\},
$$
and the entropy spectrum $\e_{\mu}(\alpha):=h_{\top}(T,K^{(\mu)}_{\alpha})$.

In \cite{Takens-Verbitskiy2000}, the authors introduced the notion of \textit{weak entropy doubling condition}
$$
C_n(\epsilon)=\sup_{x}\frac{\mu(B_n(x,\epsilon))}{\mu(B_n(x,\epsilon/2))}<\infty,\text{ and } \lim_{n}\frac 1n \log C_n(\epsilon)=0.
$$
Note that the weak entropy doubling condition implies the weak {Bowen condition} for observables $X_n:\Omega\to\R$ given by
$$
X_n(x) =-\log \mu( B_n(x,\epsilon)).
$$
Indeed, it is easy to show that
$v_{n,\epsilon}(X_n)\leq \log C_n(2\epsilon)$. 


In \cite{Takens-Verbitskiy2000} it was shown that if $\mu$ is fully supported measures
satisfying the (weak) entropy doubling condition, then for any $\alpha$ one has:
$$
\e_{\mu}(\alpha):=h_{\top}(T,K^{(\mu)}_{\alpha})\le H^*_\mu(\alpha),
$$
where $H^*_\mu(\alpha)$ is the Legendre transform of the correlation entropy function
$$
H_\mu(q)=\lim_{\epsilon\to 0}\liminf_{n\to\infty}-\frac 1n \log I_\mu(q,n,\epsilon),\quad
I_\mu(q,n,\epsilon)=\begin{cases} \int\mu(B_n(x,\epsilon))^{q-1}d\mu, & \text{if } q\ne 1,\\
\int\log \mu(B_n(x,\epsilon)) d\mu,&\text{if } q=1.\\
\end{cases}
$$
Note the similarity between this result and the upper bound (Theorem A) --
under very similar mild regularity assumptions one obtains an upper bound on the multifractal spectrum.

In order to obtain a lower bound on the multifractal spectrum of local entropies
one needs additional assumptions on the measure $\mu$. For example, $\mu$ is the so-called Bowen-Gibbs measure:
$$
\frac 1{C(\epsilon)}\le \frac{\mu( B_n(x,\epsilon))}{\exp(\sum_{k=0}^{n-1} \phi(T^kx))}\le C(\epsilon).
$$
Then the multifractal analysis of level sets of local entropies
is reduced to the analysis of level sets of Birkhoff's averages 
$K_\alpha =\{\omega: \frac 1n\sum_{k=0}^{n-1} \phi(T^k\omega)\to \alpha\}$,
see \cite{Takens-Verbitskiy2000}.

\subsection{Weighted ergodic averages.}

Fan (\cite{Fan2021}) recently studied  multifractal spectra of a sequence $\{X_n\}$ in the form $X_n:=\sum_{i=0}^{n-1}f_i\circ T^{i}$, where $f_i,\;\; i\in \Z_+$ are continuous functions on the shift space $\Omega:=\A^{\Z_+}$, where $\A$ is a finite alphabet with $\#\A = q$. A motivation to study this kind of sequence has been the weighted ergodic averages $S^{w}_nf=\sum_{i=0}^{n-1}w_i f\circ T^i,\; n\in\Z_+$ of a potential $f\in C(\Omega)$, where $\{w_i:i\in\Z_+\}\subset \R$ are some fixed weights. This is indeed a special case of the former sequence if one sets $f_i=w_i f$ for all $i\in\Z_+$. To find the dimension spectra of the sequence $\{X_n\}$ in terms of the partition function, Fan made two assumptions:
\begin{align*}
&\text{(H1):} \;\;\; \text{ The limit }  \phi(t)=\lim_{n\to\infty} \frac{1}{n}\log \int_\Omega \exp{(tX_n)}d\mu \text{ exists for all } t\in\R, \\
&\text{(H2):} \;\;\; \sup_{n\in\Z_+} \sup_{x_{0}^{n-1}=y_0^{n-1}} \sum_{i=0}^{n-1} | f_i(T^i(x))-f_i(T^i(y)) |<+\infty.
\end{align*}
Recall that in the shift space, for all $Z\subset \Omega$, $h_{\top}(T, Z)=\log q \cdot \dim_H(Z)$.
Under the assumptions (H1)-(H2), Fan obtained the following theorem reformulated in terms of topological entropies:
\begin{theorem}[\cite{Fan2021} ] If $\lambda\geq 0$,
$$h_{\top}(T,\Omega)-\phi^*(\max\partial \phi(\lambda))\leq h_{\top}(T,K_{\partial \phi(\lambda)})\leq h_{\top}(T,\Omega)-\phi^*(\min\partial \phi(\lambda)),$$
where $\partial\phi(\lambda)$ is the set of all subgradients of $\phi$ at $\lambda$, and
$$
K_{[a;b]}=\Big\{\omega\in\Omega: a\leq \liminf_{n\to\infty}\frac{X_n(\omega)}{n}\leq\limsup_{n\to\infty}\frac{X_n(\omega)}{n}\leq b \Big\}\text{ for } \; a,b\in\R.
$$
If $\lambda<0$, we have similar estimates but we have to exchange the roles of $\min \partial \phi(\lambda)$ and $\max\partial\phi(\lambda)$. In particular, if $\phi$ is differentiable at $\lambda$, then for the level set $K_{\alpha}$, $\alpha=\phi'(\lambda)$, one has
\begin{equation}\label{Aihua's result}
h_{\top}(T,K_{\alpha})=h_{\top}(T,\Omega)-\phi^*(\alpha).
\end{equation}
\end{theorem}
Now we should note that condition (H2) implies the Bowen condition, and thus the weak Bowen condition as well. In fact, the condition (H2) is strictly stronger than the Bowen condition as the following example shows. 
\begin{example}
Let $\A=\{0,1\}$, and $a_n=\frac{(-1)^n}{n^{\gamma}}, \; \gamma\in (1,2], \; n\in\N$. Define $f(\omega):=\sum_{n\in\N}a_n\omega_n$ for $\omega\in \Omega$, then $f\in C(\Omega)$. It is easy to notice that 
\begin{itemize}
   \item[1)] $\sup_{n\in\N} \sum _{j=1}^{\infty}(|a_j|+...+|a_{j+n-1}|)=\sum_{j\in\N}j|a_j|=+\infty$;
   \item[2)] $\sup_{n\in\N}\sum_{j=1}^{\infty}|a_j+...+a_{j+n-1}|<+\infty$.
\end{itemize}
Thus
\begin{equation}\label{unboundedness}
\sup_{n\in\N}\sup_{\omega_0^{n-1}=\bar{\omega}_0^{n-1}}\sum_{i=0}^{n-1}|(f\circ T^i)(\omega)-(f\circ T^i)(\bar{\omega})|=+\infty,
\end{equation}
but $\sup\limits_{n\in\N} \;v_{n} (f_n)<+\infty$, where $f_n:=\sum_{i=0}^{n-1}f\circ T^i$.
\end{example}
Since the sequence $\{X_n\}$ satisfies the weak Bowen condition, Theorem B generalizes the inequality $ h_{\top}(T, K_{\alpha})\leq h_{\top}(T,\Omega)-\phi^*(\alpha)$ which is a part of (\ref{Aihua's result}) since it relaxes both assumptions (H1) and (H2). 

It should be stressed that Theorem C is not applicable in this setup since the sequence $\{X_n\}_{n\in\N}$ is not almost additive in general.

 


\section{Final remarks}
A) The results of this paper can be extended to the case that the reference measure is a Bowen-Gibbs measure, i.e., there exists potential $\zeta\in C(\Omega,\R)$ such that for every $\epsilon>0$ there exists a constant $C(\epsilon)\geq 1$ satisfying for all $n\in\N$ and $x\in\Omega$ that
\begin{equation}\label{Bowen-Gibbs property}
    C(\epsilon)^{-1} \exp( \zeta_n(x))
\le \nu( B_n(x,\epsilon))\le
C(\epsilon) \exp( \zeta_n(x)).
\end{equation}

It is easy to see that the Ahlfors-Bowen measures satisfy the above property (\ref{Bowen-Gibbs property}) with a constant potential $\zeta\equiv const$. To extend the results for such measures, instead of the topological entropy, one can consider another dimension characteristics of the sets related to the potential $\zeta$. We note that the Pesin-Pitskel topological pressure \cite{Pesin-Pitskel1984} associated with the potential $\zeta$ is suitable for this purpose. 

B) It would be interesting if one could obtain an analogue of the above results for the Hausdorff dimension instead of the topological entropy. However, the Hausdorff dimension might be more difficult to treat as it is not related to the dynamics of the underlying transformation.


    



 \begin{center}
 \textbf{Acknowledgements.}
 \end{center}
The authors are grateful to Aihua Fan for his valuable discussions.
Qian Xiao is particularly grateful to the China Scholarship Council and Leiden University for their support, which made her visit to Leiden University possible.

\end{document}